\documentclass{article}

\usepackage{microtype}
\usepackage{graphicx}
\usepackage{booktabs} % for professional tables

\usepackage{hyperref}

\usepackage[accepted]{icml2024}

\usepackage{amsmath}
\usepackage{amssymb}
\usepackage{mathtools}
\usepackage{amsthm}

\usepackage{amsfonts}
\usepackage{multirow}
\usepackage{multicol}

\usepackage{color}
\usepackage{algorithm}
\usepackage{algorithmic}

\usepackage{xcolor,colortbl}
\definecolor{LightCyan}{rgb}{0.88,1,1}

\usepackage{subfig}
\usepackage{tcolorbox}

\usepackage[capitalize,noabbrev]{cleveref}

\theoremstyle{plain}
\newtheorem{theorem}{Theorem}[section]

\newtheorem{lemma}[theorem]{Lemma}

\theoremstyle{definition}
\newtheorem{definition}[theorem]{Definition}
\newtheorem{assumption}[theorem]{Assumption}
\theoremstyle{remark}
\newtheorem{remark}[theorem]{Remark}

\usepackage[textsize=tiny]{todonotes}

\icmltitlerunning{Optimal Hessian/Jacobian-Free Nonconvex-PL Bilevel Optimization}

\begin{document}

\twocolumn[
\icmltitle{Optimal Hessian/Jacobian-Free Nonconvex-PL Bilevel Optimization}



\icmlsetsymbol{equal}{*}

\begin{icmlauthorlist}
\icmlauthor{Feihu Huang}{1,2}
%\icmlauthor{Firstname1 Lastname1}{equal,yyy}
%\icmlauthor{Firstname2 Lastname2}{equal,yyy,comp}
%\icmlauthor{Firstname3 Lastname3}{comp}
%\icmlauthor{Firstname4 Lastname4}{sch}
%\icmlauthor{Firstname5 Lastname5}{yyy}
%\icmlauthor{Firstname6 Lastname6}{sch,yyy,comp}
%\icmlauthor{Firstname7 Lastname7}{comp}
%%\icmlauthor{}{sch}
%\icmlauthor{Firstname8 Lastname8}{sch}
%\icmlauthor{Firstname8 Lastname8}{yyy,comp}
\end{icmlauthorlist}

\icmlaffiliation{1}{College of Computer Science and Technology, Nanjing University of Aeronautics and Astronautics, Nanjing, China}
\icmlaffiliation{2}{MIIT Key Laboratory of Pattern Analysis and Machine Intelligence, Nanjing, China}

\icmlcorrespondingauthor{Feihu Huang}{huangfeihu2018@gmail.com}

\icmlkeywords{Bilevel Optimization, Nonconvex, Hessian-Free, Optimal Gradient Complexity }

\vskip 0.3in
]



\printAffiliationsAndNotice{}  % leave blank if no need to mention equal contribution

\begin{abstract}
Bilevel optimization is widely applied in many machine learning tasks such as hyper-parameter learning, meta learning and reinforcement learning. Although many algorithms recently have been developed to solve the bilevel optimization problems, they generally rely on the (strongly) convex lower-level problems. More recently, some methods have been proposed to solve the nonconvex-PL bilevel optimization problems, where their upper-level problems are possibly nonconvex, and their lower-level problems are also possibly nonconvex while satisfying  Polyak-{\L}ojasiewicz (PL) condition. However, these methods still have a high convergence complexity or a high computation complexity such as requiring compute expensive Hessian/Jacobian matrices and its inverses.
In the paper, thus, we propose an efficient Hessian/Jacobian-free method (i.e., HJFBiO) with the optimal convergence complexity to solve the nonconvex-PL bilevel problems.
Theoretically, under some mild conditions, we prove that our HJFBiO method obtains an optimal convergence rate of $O(\frac{1}{T})$, where $T$ denotes the number of iterations, and has an optimal gradient complexity of $O(\epsilon^{-1})$ in finding an $\epsilon$-stationary solution. We conduct some numerical experiments on the bilevel PL game and hyper-representation learning task to demonstrate efficiency of our proposed method.
\end{abstract}

\section{Introduction}
Bilevel optimization~\citep{colson2007overview,liu2021investigating}, as an effective two-level hierarchical optimization paradigm,
is widely applied in many machine learning tasks such as hyper-parameter learning~\citep{franceschi2018bilevel},
 meta learning~\citep{franceschi2018bilevel,ji2021bilevel} and reinforcement learning~\citep{hong2020two,chakraborty2023aligning}.
In the paper, we consider a class of nonconvex bilevel optimization problems:
\begin{align}
& \min_{x \in \mathbb{R}^d, y\in y^*(x)} \ f(x,y) + \phi(x)  & \mbox{(Upper-Level)} \label{eq:1} \\
& \mbox{s.t.} \ y^*(x) \equiv \arg\min_{y\in \mathbb{R}^p} \ g(x,y),  & \mbox{(Lower-Level)} \nonumber
\end{align}
where the upper-level function $f(x,y)$ with $y\in y^*(x)$ is possibly nonconvex, and $\phi(x)$ is a convex but possibly nonsmooth regularization such as $\phi(x)=0$ when $x \in \mathcal{X}\subseteq \mathbb{R}^d$ with convex set $\mathcal{X}$ otherwise $\phi(x)=+\infty$, or $\phi(x)=\|x\|_1$. The lower-level function $g(x,y)$ is possibly nonconvex on any $y$ and satisfies  Polyak-{\L}ojasiewicz (PL) condition~\citep{polyak1963gradient}, which
 relaxes the strong convexity. The PL condition is widely used to some
  machine learning models such as the over-parameterized deep neural networks~\citep{frei2021proxy,song2021subquadratic}.
In fact, Problem~(\ref{eq:1}) widely appears in many machine learning tasks
such as meta learning~\citep{huang2023momentum} and reinforcement learning~\citep{chakraborty2023aligning}.

\begin{table*}
  \centering
  \caption{ Comparison of \textbf{gradient (or iteration) complexity }
between our method and the existing methods in solving bilevel problem~(\ref{eq:1}) for finding an $\epsilon$-stationary solution ($\|\nabla F(x)\|^2\leq \epsilon$ or its equivalent variants, where $F(x)=f(x,y)$ with $y\in y^*(x)$). Here $g(x,\cdot)$ denotes function on the second variable $y$ with fixing variable $x$.  \textbf{SC} stands for strongly convex. \textbf{H.J.F.} stands for Hessian/Jacobian-Free. \textbf{L.H.} stands for Lipschitz Hessian condition. \emph{The Prox-F$^2$BA and F$^2$BA methods rely on some strict conditions such as Lipschitz Hessian of function $f(x,y)$.} \textbf{Note that} the GALET~\citep{xiao2023generalized} method simultaneously uses the PL condition, its Assumption 2 (i.e., let $\sigma_g = \inf_{x,y}\{\sigma_{\min}^{+}(\nabla^2_{yy} g(x,y))\} >0$ for all $(x,y)$) and its Assumption 1 (i.e., $\nabla^2_{yy} g(x,y)$ is Lipschitz continuous). Clearly, when Hessian matrix $\nabla^2_{yy} g(x,y)$ is singular, its Assumption 1 and Assumption 2 imply that the lower bound of the non-zero singular values $\sigma_g$ is close to zero (i.e., $\sigma_g\rightarrow 0$), under this case, the convergence results of the GALET are \textbf{meaningless}, e.g., the constant $L_w = \frac{\ell_{f,1}}{\sigma_g}+\frac{\sqrt{2}\ell_{g,2}\ell_{f,0}}{\sigma_g^2}\rightarrow + \infty$ used in its Lemmas 6 and 9. Under the other case, the PL condition, Lipschitz continuous of Hessian and its Assumption 2 (the singular values of Hessian is bounded away from 0, i.e., $\sigma_g>0$) imply that GALET assumes strongly convex (Detailed discussion in the Appendix~\ref{appendix-B}).   }
  \label{tab:1}
 \vspace*{-4pt}
  % \resizebox{0.98\textwidth}{!}{
\begin{tabular}{c|c|c|c|c|c|c}
  \hline
  \textbf{Algorithm} & \textbf{Reference} &  $g(x,\cdot)$ &  \textbf{L.H.} on $f(\cdot,\cdot)$ & \textbf{Complexity} & \textbf{Loop(s)}  & \textbf{H.J.F.}  \\ \hline
   BOME  & \cite{liu2022bome} & PL / local-PL &  & $O(\epsilon^{-1.5})$ / $O(\epsilon^{-2})$ & Double & $\surd$  \\   \hline
  V-PBGD  & \cite{shen2023penalty} & PL / local-PL &  & $O(\epsilon^{-1.5})$ / $O(\epsilon^{-1.5})$ & Double & $\surd$  \\  \hline
   GALET  & \cite{xiao2023generalized} & SC / PL &  & $O(\epsilon^{-1})$ / \textbf{Meaningless} & Triple &  \\ \hline
   SLM & \cite{lu2023slm} & PL / local-PL &  & $O(\epsilon^{-3.5})$ / $O(\epsilon^{-3.5})$ & Double  & $\surd$  \\ \hline
  Prox-F$^2$BA  & \cite{kwon2023penalty} & Proximal-EB & $\surd$ &  $O(\epsilon^{-1.5})$ / $O(\epsilon^{-1.5})$ & Double & $\surd$ \\ \hline
  F$^2$BA  & \cite{chen2024bilevel} & PL / local-PL & $\surd$ &  $O(\epsilon^{-1})$ / $O(\epsilon^{-1})$ & Double &  $\surd$ \\ \hline
  MGBiO  & \cite{huang2023momentum} & PL / local-PL &   &  $O(\epsilon^{-1})$ / $O(\epsilon^{-1})$ & Single &  \\ \hline
  AdaPAG  & \cite{huang2023adaptive} & PL / local-PL &  &  $O(\epsilon^{-1})$ / $O(\epsilon^{-1})$ & Single &  \\ \hline
HJFBiO  & Ours  & PL / local-PL & & \color{red}{$O(\epsilon^{-1})$} / \color{red}{$O(\epsilon^{-1})$} & \color{red}{Single} & \color{red}{$\surd$}  \\ \hline
\end{tabular}
% }
 \vspace*{-6pt}
\end{table*}

The inherent nested nature of bilevel optimization gives rise to several difficulties
in effectively solving these bilevel problems.
For example, compared with the standard single-level optimization (i.e., $g\big(x,y\big)=0$ in Problem~(\ref{eq:1})), the main difficulty of bilevel optimization
is that the minimization of the upper and lower-level objectives are intertwined via the
minimizer $y^*(x)\in \arg\min_y g(x,y)$ of the lower-level problem.
To deal with this difficulty, recently many bilevel optimization methods~\citep{ghadimi2018approximation,hong2020two,ji2021bilevel,huang2022enhanced,chen2023accelerated} have been proposed
by imposing the
strong convexity assumption on the Lower-Level (LL) problems.
 The LL strong convexity assumption ensures the uniqueness of LL minimizer (i.e., LL Singleton),
which simplifies both the optimization process and theoretical analysis, e.g., hyper-gradient $\nabla F(x)$ of
 the upper-level objective $F(x)=f(x,y^*(x))$ has a simple closed-form:
\begin{align} \label{eq:2}
 \nabla F(x) & = \nabla_xf(x,y^*(x)) \\
  & - \nabla^2_{xy}g(x,y^*(x))\nabla^2_{yy}g(x,y^*(x))^{-1}\nabla_yf(x,y^*(x)). \nonumber
\end{align}
Based on the above form of hyper-gradient $\nabla F(x)$, some gradient-based methods~\citep{chen2022single,dagreou2022framework} require an explicit extraction of second-order information of $g(x,y)$ with a major focus on efficiently estimating its Jacobian and inverse Hessian.
Meanwhile, the other gradient-based methods~\citep{li2022fully,dagreou2022framework,sow2022convergence,yang2023achieving} avoid directly estimating its second-order computation and only use the first-order information of both upper and lower objectives.

Recently, to relax the LL strong convexity assumption, another line of research
is dedicated to bilevel optimization with convex LL problems, which bring about several challenges such as the
presence of multiple LL local optimal solutions (i.e., Non-Singleton). Under this case, there does not exist the above  hyper-gradient form~(\ref{eq:2}). To handle this concern, some effective methods~\citep{sow2022constrained,liu2023averaged,lu2023first,cao2023projection} recently have been developed. For example, \cite{sow2022constrained} developed the primal-dual algorithms for bilevel optimization with multiple inner minima in the LL problem. Subsequently, \cite{lu2023first} studied the constrained bilevel optimization with convex lower-level. Meanwhile, \cite{liu2023averaged} proposed an effective averaged method of multipliers for bilevel optimization with convex lower-level.

In fact, the bilevel optimization problems with nonconvex LL problems frequently appear in many machine tasks such as hyper-parameter learning in training deep neural networks.
Since the above bilevel optimization methods mainly rely on the
restrictive LL strong convexity or convexity assumption, clearly, they can not effectively solve the
bilevel optimization problems with nonconvex LL problems. Recently, some bilevel approaches~\citep{liu2021value,liu2022bome,
chen2023bilevel,liu2023averaged,huang2023momentum,huang2023adaptive,kwon2023penalty,chen2024bilevel} studied the bilevel optimization with non-convex lower-level.
For example,
\cite{liu2022bome} proposed an effective first-order method for nonconvex-PL bilevel optimization, where the lower-level problem is nonconex but satisfies PL condition. \cite{shen2023penalty} designed an effective penalty-based gradient method for the constrained nonconvex-PL bilevel optimization. \cite{kwon2023penalty} studied the nonconvex bilevel optimization with nonconvex lower-level satisfying  proximal error-bound (EB) condition that is analogous to PL condition.
Meanwhile, \cite{huang2023momentum} proposed a class of efficient momentum-based gradient methods for the nonconvex-PL bilevel optimization, which obtain an optimal gradient complexity but rely on
requiring compute expensive projected Hessian/Jacobian matrices and its inverses. Subsequently, \cite{xiao2023generalized} proposed a generalized alternating method (i.e., GALET)  for nonconvex-PL bilevel optimization, which still relies on the expensive  Hessian/Jacobian matrices. Unfortunately, the convergence results of the GALET method are meaningless (please see Table~\ref{tab:1}).
 Thus, there exists a natural question:
\begin{center}
\begin{tcolorbox}
\textbf{ Could we propose an efficient Hessian/Jacobian-free method for the nonconvex-PL bilevel optimization with an optimal  complexity? }
\end{tcolorbox}
\end{center}
In the paper, we affirmatively answer to this question, and propose
an efficient Hessian/Jacobian-free method
to solve Problem~(\ref{eq:1}) based on the finite-difference estimator and
a new projection operator.
Our main contributions are given:
\begin{itemize}
\item[(i)] We propose an efficient Hessian/Jacobian-free method (i.e., HJFBiO)  based on the finite-difference estimator and a new useful projection operator. In particular, our HJFBiO method not only uses low computational first-order gradients instead of high computational Hessian/Jacobian matrices, but also  applies the low computational new projection operator to
    vector variables instead of matrix variables. Thus, our HJFBiO method has a lower computation at each iteration.
\item[(ii)] We provide a solid convergence analysis for our HJFBiO method. Under some mild conditions, we prove that our HJFBiO method reaches the best known iteration (gradient) complexity of $O(\epsilon^{-1})$ for finding an $\epsilon$-stationary solution of Problem~(\ref{eq:1}), which matches the lower bound established by the first-order method for finding an $\epsilon$-stationary point of nonconvex smooth optimization problems~\citep{carmon2020lower}.
\item[(iii)] We conduct some numerical experiments including bilevel Polyak-{\L}ojasiewicz game and hyper-representation learning to demonstrate efficiency of our proposed method.
\end{itemize}
Meanwhile,  \cite{chen2024bilevel} proposed a F$^2$BA method for the nonconvex-PL bilevel optimization, which also obtains the best known gradient complexity $O(\epsilon^{-1})$ for finding an $\epsilon$-stationary solution of Problem~(\ref{eq:1}), but it relies on some stricter conditions such as Lipschitz Hessian of the upper function $f(x,y)$. Under these strict conditions, although the F$^2$BA method~\cite{chen2024bilevel} obtains a gradient complexity $O(\epsilon^{-1})$,
 this is not an optimal gradient complexity (Detailed discussion in the Appendix~\ref{appendix-B}).

\subsection*{Notations}
Given function $f(x,y)$, $f(x,\cdot)$ denotes  function \emph{w.r.t.} the second variable with fixing $x$,
and $f(\cdot,y)$ denotes function \emph{w.r.t.} the first variable
with fixing $y$. $\nabla_x$ denotes the partial derivative on variable $x$. Let $\nabla^2_{xy}=\nabla_x\nabla_y$ and $\nabla^2_{yy}=\nabla_y\nabla_y$. $\|\cdot\|$ denotes the $\ell_2$ norm for vectors and spectral norm for matrices. $\langle x,y\rangle$ denotes the inner product of two vectors $x$ and $y$. $I_{d}$ denotes a $d$-dimensional identity matrix. $a_t=O(b_t)$ denotes that $a_t \leq c b_t$ for some constant $c>0$. The notation $\tilde{O}(\cdot)$ hides logarithmic terms.

$\mathcal{S}_{[\mu,L_g]}[\cdot]$ denotes a
projection on the set $\{X\in \mathbb{R}^{d\times d}: \mu \leq \varrho(X) \leq L_g\}$,
where $\varrho(\cdot)$ denotes the eigenvalue function. $\mathcal{S}_{[\mu,L_g]}$ can be implemented by using Singular Value Decomposition (SVD) and thresholding the singular values.
$\mathcal{P}_{r_v}(\cdot)$ is a projection onto set $\{v\in \mathbb{R}^p: \|v\| \leq r_v>0\}$.
\section{ Preliminaries }
In this section, we provide some mild assumptions and useful lemmas on the above Problem~(\ref{eq:1}).
\subsection{ Mild Assumptions}
\begin{assumption} \label{ass:1}
 The function $g(x,\cdot)$ satisfies the Polyak-{\L}ojasiewicz (PL) condition,
 if there exist $\mu > 0$ such that for any given $x$, it holds that
\begin{align}
 \|\nabla_y g(x,y)\|^2 \geq 2\mu \big(g(x,y)-\min_y g(x,y)\big), \ \forall y\in \mathbb{R}^p. \nonumber
\end{align}
\end{assumption}

\begin{assumption} \label{ass:2}
 The function $g(x,y)$ is nonconvex and satisfies
\begin{align}
 {\color{blue}{\varrho\big(\nabla^2_{yy}g\big(x,y^*(x)\big)\big) \in [\mu, L_g]}},
\end{align}
where $y^*(x) \in \arg\min_y g(x,y)$, and $\varrho(\cdot)$ denotes the eigenvalue (or singular-value) function and $L_g\geq\mu>0$.
\end{assumption}

\begin{assumption} \label{ass:3}
The functions $f(x,y)$ and $g(x,y)$ satisfy
\begin{itemize}
 \item[1)] For all $x,y$, we have $ \|\nabla_yf(x,y)\|\leq C_{fy}, \ \|\nabla^2_{xy}g(x,y)\|\leq C_{gxy}$;
 \item[2)] The partial derivatives $\nabla_x f(x,y)$ and $\nabla_y f(x,y)$ are $L_f$-Lipschitz continuous;
 \item[3)] The partial derivatives $\nabla_xg(x,y)$ and $\nabla_yg(x,y)$ are $L_g$-Lipschitz continuous.
\end{itemize}
\end{assumption}

\begin{assumption} \label{ass:4}
 The partial derivatives $\nabla^2_{xy}g(x,y)$ and $\nabla^2_{yy}g(x,y)$ are $L_{gxy}$-Lipschitz and $L_{gyy}$-Lipschitz, e.g.,
 for all $x,x_1,x_2 \in \mathbb{R}^d$ and $y,y_1,y_2 \in \mathbb{R}^p$
 \begin{align}
   &\|\nabla^2_{xy} g(x_1,y)-\nabla^2_{xy} g(x_2,y)\| \leq L_{gxy}\|x_1-x_2\|,  \nonumber \\
   & \|\nabla^2_{xy} g(x,y_1)-\nabla^2_{xy} g(x,y_2)\| \leq L_{gxy}\|y_1-y_2\|. \nonumber
 \end{align}
 \end{assumption}

\begin{assumption} \label{ass:5}
 The function $\Phi(x)=F(x)+\phi(x)$ is bounded below in $x\in \mathbb{R}^d$, \emph{i.e.,} $\Phi^* = \inf_{x\in \mathbb{R}^d}\Phi(x) > -\infty$.
\end{assumption}
Assumption~\ref{ass:1} is commonly used in bilevel optimization without the lower-level strongly convexity~\citep{liu2022bome,shen2023penalty,huang2023momentum}.
 Assumption~\ref{ass:2} ensures that the minimizer $y^*(x)=\arg\min_y g(x,y)$
 is unique,
 which imposes the non-singularity of $\nabla^2_{yy}g(x,y)$ only at the minimizers $y^*(x) \in \arg\min_y g(x,y)$, as in~\citep{huang2023momentum}. \emph{Based on Lemma G.6 of \cite{chen2024bilevel} given in the Appendix~\ref{appendix-B},
 our Assumption~\ref{ass:2} is reasonable when
 has an unique minimizer.} Meanwhile, we also study the case that $\min_y g(x,y)$
 has multiple local minimizers in the following section~\ref{sec4.2}. \textbf{Note that} since $y^*(x)\in \arg\min_y g(x,y)$, we can not have negative eigenvalues at the minimizer $y^*(x)$, so Assumption~\ref{ass:2}  assumes that $\varrho\big(\nabla^2_{yy}g\big(x,y^*(x)\big)\big) \in [\mu, L_g]$ instead of $\varrho\big(\nabla^2_{yy}g\big(x,y^*(x)\big)\big) \in [-L_g,-\mu] \cup [\mu, L_g]$.
Since $\nabla^2_{yy} g(x,y)$
is a symmetric matrix, its singular values are the absolute value of eigenvalues.
Hence, we also can use $\varrho(\cdot)$ to denote the singular-value function.

Assumption~\ref{ass:3} is commonly appeared in bilevel optimization methods \citep{ghadimi2018approximation,ji2021bilevel,liu2022bome}.
 Meanwhile, the BOME~\cite{liu2022bome} uses the stricter assumption that
  $\|\nabla f(x,y)\|$, $\|\nabla g(x,y)\|$, $|f(x,y)|$ and $|g(x,y)|$
 are bounded for any $(x,y)$ in its Assumption~3. Assumption~\ref{ass:4} is also commonly used in bilevel optimization methods \citep{ghadimi2018approximation,ji2021bilevel}.
Assumption~\ref{ass:5} ensures the feasibility of the bilevel Problem~(\ref{eq:1}).

 For example, we consider a nonconvex-PL  bilevel problem
\begin{align}
\min_{x\in [1,2],y\in y^*(x)} & \Big\{f(x,y)=x^2 + y^2 + 3x\sin^2(y)\Big\}, \\
\mbox{s.t.} \ & y^*(x) \equiv\min_{y\in \mathbb{R}} \Big\{g(x,y)=xy^2+x\sin^2(y)\Big\},  \nonumber
\end{align}
which can be rewritten as
\begin{align}
\min_{x\in \mathbb{R},y\in y^*(x)} & \Big\{f(x,y)+ \phi(x)\Big\},  \label{eq:exp} \\
\mbox{s.t.} & \ y^*(x) \equiv\min_{y\in \mathbb{R}} \Big\{g(x,y)=xy^2+x\sin^2(y)\Big\}, \nonumber
\end{align}
where $\phi(x)=0$ when $x\in [1,2]$ otherwise $\phi(x)=+\infty$.
From the above bilevel problem~(\ref{eq:exp}), we can easily obtain $y^*(x)=0$, $\nabla^2_{yy}g(x,y)=x\big(2+2\cos^2(y)-2\sin^2(y)\big)$, and then we have $\nabla^2_{yy}g(x,y^*(x))=4x>0$ due to $x\in [1,2]$. Since $\nabla^2_{yy}g(x,y^*(x))=4x>0$, our Assumption~\ref{ass:2} holds. Meanwhile  $\nabla^2_{yy}g(x,y)=x\big(2+2\cos^2(y)-2\sin^2(y)\big)$ for any $y\in \mathbb{R}$ may be zero or negative such as $\nabla^2_{yy} g(x,\pi/2)=0$. Meanwhile, for any $y\in \mathbb{R}$, clearly $|f(x,y)|$ and $|g(x,y)|$ are not bounded. Thus, the assumptions used in~\citep{liu2022bome,xiao2023generalized,kwon2023penalty} may be not satisfied.

\subsection{Useful Lemmas}
In this subsection, based on the above assumptions, we give some useful lemmas.

\begin{lemma} \label{lem:1}
(\cite{huang2023momentum})
Under the above Assumption \ref{ass:2}, we have, for any $x\in \mathbb{R}^d$,
\begin{align}
 &\nabla F(x) =\nabla_x f(x,y^*(x)) \nonumber \\
 &  \qquad - \nabla^2_{xy} g(x,y^*(x))\Big[\nabla^2_{yy}g(x,y^*(x))\Big]^{-1}\nabla_y f(x,y^*(x)). \nonumber
\end{align}
\end{lemma}

From the above Lemma~\ref{lem:1}, we can get the same form of hyper-gradient $\nabla F(x)$ as in~(\ref{eq:2}). Since the Hessian matrix $\nabla^2_{yy}g(x,y)$ for all $(x,y)$ may be singular, as in \cite{huang2023momentum},
we define a useful hyper-gradient estimator:
\begin{align}
& \hat{\nabla} f(x,y)=\nabla_xf(x,y) \nonumber \\
& \qquad \quad - \nabla^2_{xy}g(x,y) \big({\color{blue}{\mathcal{S}_{[\mu,L_g]}\big[\nabla^2_{yy}g(x,y)\big]}}\big)^{-1} \nonumber
\nabla_yf(x,y),
\end{align}
which replaces the standard hyper-gradient estimator $\breve{\nabla} f(x,y)$ used in ~\cite{ghadimi2018approximation,ji2021bilevel} for the strongly-convex lower-level optimization,
\begin{align}
\breve{\nabla} f(x,y) & =\nabla_xf(x,y) \nonumber \\
& \quad  \quad - \nabla^2_{xy}g(x,y) \big(\nabla^2_{yy}g(x,y)\big)^{-1}\nabla_yf(x,y). \nonumber
\end{align}

\begin{lemma} \label{lem:2}
(\cite{huang2023momentum})
Under the above Assumptions \ref{ass:1}-\ref{ass:4}, the functions (or mappings) $F(x)=f(x,y^*(x))$, $G(x)=g(x,y^*(x))$ and $y^*(x)\in \arg\min_{y\in \mathbb{R}^p}g(x,y)$ satisfy, for all $x_1,x_2\in \mathbb{R}^d$,
\begin{align}
 & \|y^*(x_1)-y^*(x_2)\| \leq \kappa\|x_1-x_2\|, \nonumber \\
 & \|\nabla y^*(x_1) - \nabla y^*(x_2)\| \leq L_y\|x_1-x_2\|, \nonumber \\
 & \|\nabla F(x_1) - \nabla F(x_2)\|\leq L_F\|x_1-x_2\|, \nonumber \\
 & \|\nabla G(x_1) - \nabla G(x_2)\|\leq L_G\|x_1-x_2\|, \nonumber
\end{align}
where $\kappa=C_{gxy}/\mu$, $L_y=\big( \frac{C_{gxy}L_{gyy}}{\mu^2} +  \frac{L_{gxy}}{\mu} \big) (1+ \frac{C_{gxy}}{\mu})$,
$L_F=\Big(L_f + L_f\kappa + C_{fy}\big( \frac{C_{gxy}L_{gyy}}{\mu^2} +  \frac{L_{gxy}}{\mu} \big)\Big)(1+\kappa)$ and $L_G=\Big(L_g + L_g\kappa + C_{gy}\big( \frac{C_{gxy}L_{gyy}}{\mu^2} +  \frac{L_{gxy}}{\mu} \big)\Big)(1+\kappa)$.
\end{lemma}

\begin{lemma} \label{lem:3}
(\cite{huang2023momentum})
Let $\nabla F(x) = \nabla f(x,y^*(x))$ and $\hat{\nabla} f(x,y)=\nabla_xf(x,y) - \nabla^2_{xy}g(x,y) \big({\color{blue}{\mathcal{S}_{[\mu,L_g]}\big[\nabla^2_{yy}g(x,y)\big]}}\big)^{-1}\nabla_yf(x,y)$,
 we have
 \begin{align}
 \|\hat{\nabla} f(x,y)-\nabla F(x)\|^2  \leq \frac{2\hat{L}^2}{\mu}\big(g(x,y)-\min_y g(x,y)\big), \nonumber
\end{align}
where $\hat{L}^2 = 4\big(L^2_f+ \frac{L^2_{gxy}C^2_{fy}}{\mu^2} + \frac{L^2_{gyy} C^2_{gxy}C^2_{fy}}{\mu^4} +
 \frac{L^2_fC^2_{gxy}}{\mu^2}\big)$.
\end{lemma}
\section{ Efficient Hessian/Jacobian-Free Bilevel Optimization Method }
In the section, we  propose an efficient Hessian/Jacobian-free method
to solve the nonconvex-PL bilevel Problem~(\ref{eq:1})
based on the finite-difference estimator and
a new projection operator.
Here we first define a useful projection operator:

\begin{definition}\label{def:1}
 Given matrix $H\in \mathbb{R}^{p\times p}$ and vector $v\in \mathbb{R}^p$, and $\mathcal{S}_{[\mu,L_g]}[\cdot]$
 is a projection operator on the set $\{H\in \mathbb{R}^{p\times p}: 0<\mu \leq \varrho(H) \leq L_g\}$ where $\varrho(\cdot)$ denotes the eigenvalue function, and $\mathcal{P}_{r_v}(\cdot)$ is a projection operator
 onto the set $\{v\in \mathbb{R}^p: \|v\| \leq r_v>0\}$, then we define a \textbf{new projection operator} $\mathcal{M}_{r_h}(\cdot,\cdot)$ on set $\{H\in \mathbb{R}^{p\times p}, v\in \mathbb{R}^p: \|Hv\| \leq r_h\}$, which satisfies
 \begin{align}
  \mathcal{M}_{r_h}\big(H,v\big) \triangleq \mathcal{S}_{[\mu,L_g]}[H]\mathcal{P}_{r_v}(v),
 \end{align}
 where $0<r_h\leq r_v L_g$.
\end{definition}
For notational simplicity, let $\mathcal{M}_{r_h}\big(H,v\big)=\mathcal{M}_{r_h}\big(Hv\big)$ in the following.

\begin{algorithm}[tb]
\caption{ Hessian/Jacobian-free Bilevel Optimization (i.e, HJFBiO) Algorithm }
\label{alg:1}
\begin{algorithmic}[1] %[1] enables line numbers
\STATE {\bfseries Input:} $T$, learning rates $\lambda>0$, $\gamma>0$, $\tau>0$, and tuning parameters $\delta_\epsilon>0$, $r_v>0$, $r_h>0$,
and initial input $x_1 \in \mathbb{R}^d$, $y_1 \in \mathbb{R}^p$ and $v_1 \in \mathbb{R}^p$ ; \\
\FOR{$t = 1, 2, \ldots, T$}
\STATE  Compute $u_t = \nabla_y g(x_t,y_t)$, and update $y_{t+1}=y_t-\lambda u_t$;
\STATE  Compute $w_t = \widetilde{\nabla}f(x_t, y_t, v_t) = \nabla_xf(x_t,y_t) - \widetilde{J}(x_t, y_t, v_t, \delta_{\epsilon})$, and update $x_{t+1}=\mathbb{P}^\gamma_{\phi(\cdot)}\big(x_t, w_t\big)$;
\STATE  Compute $h_t =  \widetilde{\nabla}_v R(x_t,y_t,v_t) = \mathcal{M}_{r_h}\big(\widetilde{H}(x_t, y_t, v_t, \delta_{\epsilon})\big) - \nabla_yf(x_t,y_t)$, and update $v_{t+1} = \mathcal{P}_{r_v}\big(v_t -\tau h_t\big)$;
\ENDFOR
\STATE {\bfseries Output:} Chosen uniformly random from $\{x_t\}_{t=1}^{T}$.
\end{algorithmic}
\end{algorithm}

From the above Lemma~\ref{lem:1}, the hyper-gradient $\nabla F(x)$ takes the form of
\begin{align}\label{def:Phi}
 \nabla & F(x) = \nabla_xf(x,y^*(x))  \\
 & - \nabla_{xy}^2g(x,y^*(x))\big[\nabla_{yy}^2g(x,y^*(x))\big]^{-1}\nabla_y f(x,y^*(x)). \nonumber
\end{align}
In the above problem~(\ref{eq:1}), the lower objective function $g(x,y)$ on variable $y$ is not strongly convex, so its Hessian matrix $\nabla^2_{yy}g(x,y)$ for all $(x,y)$ may be singular. As in~\cite{huang2023momentum},
 we define a useful hypergradient estimator:
\begin{align}\label{def:hg}
\widehat{\nabla} & f(x,y)=\nabla_xf(x,y) \\
& - \nabla^2_{xy}g(x,y) \big({\color{blue}{\mathcal{S}_{[\mu,L_g]}\big[\nabla^2_{yy}g(x,y)\big]}}\big)^{-1}
\nabla_yf(x,y). \nonumber
\end{align}
Since the above hypergradient~(\ref{def:hg}) requires computing the expensive projected Hessian inverse, as in~\cite{yang2023achieving}, we define a new hypergradient surrogates as follows:
\begin{align}\label{def:hgs}
    \widehat{\nabla}f(x, y, v) = \nabla_xf(x,y) - \nabla^2_{xy}g(x,y)v,
\end{align}
where $v\in\mathbb{R}^p$ is an auxiliary vector to approximate the projected Hessian-inverse vector product $\big(\mathcal{S}_{[\mu,L_g]}\big[\nabla^2_{yy}g(x,y)\big]\big)^{-1}
\nabla_yf(x,y)$ in (\ref{def:hg}), which can be rewritten as a solution of the following linear system:
\begin{align}\label{def:v}
&v^*= \big({\color{blue}{\mathcal{S}_{[\mu,L_g]}\big[\nabla^2_{yy}g(x,y)\big]}}\big)^{-1}
\nabla_yf(x,y)  \\
& =\arg\min_{v\in\mathbb{R}^p} \Big\{ \frac{1}{2}v^T{\color{blue}{\mathcal{S}_{[\mu,L_g]}\big[\nabla^2_{yy}g(x,y)\big]}}v - v^T\nabla_yf(x,y) \Big\}.
\nonumber
\end{align}
Under this case, we can use the following new iterations to solve the nonconvex-PL bilevel problem~(\ref{eq:1}):
for $t\geq1$,
\begin{align}
 \begin{cases} & y_{t+1}=y_t -\lambda\nabla_y g(x_t,y_t), \\
 &x_{t+1} = \mathbb{P}^\gamma_{\phi(\cdot)}\big(x_t,\widehat{\nabla}f(x_t, y_t, v_t)\big), \\
 &v_{t+1} = \mathcal{P}_{r_v}\big(v_t -\tau\nabla_v R(x_t,y_t,v_t)\big),
 \end{cases}
\end{align}
where $\lambda>0$, $\gamma>0$ and $\tau>0$ are learning rates, and the proximal operator defined as: given vectors $x_t,w_t\in \mathbb{R}^d$,
\begin{align}
 & \mathbb{P}^\gamma_{\phi(\cdot)}\big(x_t,w_t\big) \nonumber \\
 & = \arg\min_{x\in \mathbb{R}^d}\big\{ \langle w_t, x\rangle
+ \frac{1}{2\gamma}\|x-x_t\|^2 + \phi(x)\big\}.
\end{align}
In particular, on updating variable $v\in \mathbb{R}^p$, we use a projection
$\mathcal{P}_{r_v}(\cdot)$ onto the set $\{v\in \mathbb{R}^p: \|v\| \leq r_v\}$ with $0<r_v\leq \frac{C_{fy}}{\mu}$
to obtain the bounded variable $v$.
Here we use the function
\begin{align}
 R(x,y,v) = \frac{1}{2}v^T\mathcal{S}_{[\mu,L_g]}\big[\nabla^2_{yy}g(x,y)\big]v - v^T\nabla_yf(x,y),\nonumber
\end{align}
and then we have $\nabla_v R(x,y,v) = \mathcal{S}_{[\mu,L_g]}\big[\nabla^2_{yy}g(x,y)\big]v - \nabla_yf(x,y)$.

In the high-dimensional setting, clearly computing Hessian matrix $\nabla^2_{yy}g(x,y)$ and Jacobian matrix
$\nabla^2_{xy}g(x,y)$ is expensive. To approximate this hypergradient efficiently, we further use
the finite-difference technique to estimate the Hessian-vector $\nabla^2_{yy}g(x,y)v$ and Jacobian-vector
$\nabla^2_{xy}g(x,y)v$ products. Specifically, given a small constant $\delta_{\epsilon}>0$, we define two
finite-difference estimators to estimate $\nabla^2_{yy}g(x,y)v$ and $\nabla^2_{xy}g(x,y)v$ respectively, defined as:
\begin{align}
& \widetilde{H}(x, y, v, \delta_{\epsilon}) = \frac{\nabla_y g(x, y+\delta_{\epsilon} v) - \nabla_y g(x, y-\delta_{\epsilon} v)}{2\delta_{\epsilon}},  \label{eq:Hfd}\\
& \widetilde{J}(x, y, v, \delta_{\epsilon}) = \frac{\nabla_x g(x, y+\delta_{\epsilon} v) - \nabla_x g(x, y-\delta_{\epsilon} v)}{2\delta_{\epsilon}}.  \label{eq:Jfd}
\end{align}
Then we can use the following Hessian/Jacobian-free iterations to solve
the nonconvex-PL bilevel problem~(\ref{eq:1}):
at $(t+1)$-th iteration,
\begin{align} \label{eq:21}
\begin{cases}
 & y_{t+1}=y_t -\lambda\nabla_y g(x_t,y_t), \\
 & x_{t+1} = \mathbb{P}^\gamma_{\phi(\cdot)}\big(x_t,\widetilde{\nabla}f(x_t, y_t, v_t)\big), \\
 & v_{t+1} = \mathcal{P}_{r_v}\big(v_t -\tau\widetilde{\nabla}_v R(x_t,y_t,v_t)\big),
 \end{cases}
\end{align}
where $0<r_v\leq \frac{C_{fy}}{\mu}$, $0<r_h\leq r_v L_g \leq \frac{C_{fy} L_g}{\mu}$, $\lambda>0$, $\gamma>0$, $\tau>0$ and
\begin{align}
 & \widetilde{\nabla}f(x_t, y_t, v_t) = \nabla_xf(x_t,y_t) - \widetilde{J}(x_t, y_t, v_t, \delta_{\epsilon}), \nonumber \\
 & \widetilde{\nabla}_v R(x_t,y_t,v_t) = \mathcal{M}_{r_h}\big(\widetilde{H}(x_t, y_t, v_t, \delta_{\epsilon})\big) - \nabla_yf(x_t,y_t). \nonumber
\end{align}
Based on the above iterations~(\ref{eq:21}),  we give a procedure framework of our HJFBiO algorithm in Algorithm \ref{alg:1}. When $\phi(x)\equiv 0$, in updating the variable $x$, we have $x_{t+1}=x_t-\gamma w_t=x_t -\gamma\widetilde{\nabla}f(x_t, y_t, v_t)$.

\textbf{Note that} in our Algorithm~\ref{alg:1}, we use two \emph{low computational}
finite-difference estimators~(\ref{eq:Hfd}) and~(\ref{eq:Jfd}) instead of computing \emph{high computational}
Hessian matrix $\nabla^2_{yy}g(x,y)\in \mathbb{R}^{d\times d}$ and Jacobian matrix $\nabla^2_{xy}g(x,y)
\in \mathbb{R}^{d \times p}$. Moreover, we also use a \emph{low computational}
projection operator $\big\{||\widetilde{H}(x, y, v, \delta_{\epsilon}) || \leq r_h, 0<r_h\leq r_v L_g \leq \frac{C_{fy} L_g}{\mu} \big\}$ on vector $\widetilde{H}(x, y, v, \delta_{\epsilon})$ instead of computing \emph{high computational} projection operator
$\mathcal{S}_{[\mu,L_g]}\big[\nabla^2_{yy} g(x,y) \big]$ used in~\citep{huang2023momentum,huang2023adaptive}. Thus, our HJFBiO algorithm only requires a low computational complexity of $O(p+d)$ at each iteration.

From Algorithm \ref{alg:1}, since $v_t = \mathcal{P}_{r_v}(v_t)$, we have
\begin{align}
 &\mathcal{M}_{r_h}\big(\widetilde{H}(x_t, y_t, v_t, \delta_{\epsilon})\big) \nonumber \\
 &   =
 \mathcal{M}_{r_h}\Big(\frac{1}{2\delta_\epsilon}\big(\nabla_y g(x_t, y_t+\delta_{\epsilon} v_t) - \nabla_y g(x_t, y_t-\delta_{\epsilon} v_t)\big)\Big)   \nonumber \\
 &  =
 \mathcal{M}_{r_h}\Big(\frac{1}{2\delta_\epsilon}\big(\nabla_y g(x_t, y_t+\delta_{\epsilon} v_t) -\nabla_y g(x_t, y_t) \nonumber \\
 & \quad + \nabla_y g(x_t, y_t)- \nabla_y g(x_t, y_t-\delta_{\epsilon} v_t)\big)\Big)   \nonumber \\
  & = \mathcal{M}_{r_h}\Big(
\frac{1}{2\delta_{\epsilon}} \int_{k=0}^1\nabla^2_{yy} g(x_t, y_t+k\delta_{\epsilon} v_t) \delta_{\epsilon}v_t dk \nonumber \\
& \quad + \frac{1}{2\delta_{\epsilon}}\int_{k=0}^1\nabla^2_{yy} g(x_t, y_t-k\delta_{\epsilon} v_t) \delta_{\epsilon}v_t dk\Big)   \nonumber \\
 & = \mathcal{M}_{r_h}\bigg(
 \Big(\frac{1}{2}\int_{k=0}^1\nabla^2_{yy} g(x_t, y_t+k\delta_{\epsilon} v_t) dk \nonumber \\
& \quad + \frac{1}{2}\int_{k=0}^1\nabla^2_{yy} g(x_t, y_t-k\delta_{\epsilon} v_t) dk \Big)v_t\bigg)   \nonumber \\
 & = \mathcal{S}_{[\mu,L_g]}\Big[\frac{1}{2}\int_{k=0}^1\nabla^2_{yy} g(x_t, y_t+k\delta_{\epsilon} v_t) dk \nonumber \\
&\quad + \frac{1}{2}\int_{k=0}^1\nabla^2_{yy} g(x_t, y_t-k\delta_{\epsilon} v_t) dk \big)\Big]v_t,
\end{align}
where the last equality holds by $v_t = \mathcal{P}_{r_v}(v_t)$ and the above definition~\ref{def:1}.
Thus, we have
\begin{align}
 & \lim_{\delta_\epsilon\rightarrow 0 }\mathcal{M}_{r_h}\big(\widetilde{H}(x_t, y_t, v_t, \delta_{\epsilon})\big) \nonumber \\
 & = \lim_{\delta_\epsilon\rightarrow 0 }\mathcal{S}_{[\mu,L_g]}\Big[\frac{1}{2}\int_{k=0}^1\nabla^2_{yy} g(x_t, y_t+k\delta_{\epsilon} v_t) dk \nonumber \\
 & \quad + \frac{1}{2}\int_{k=0}^1\nabla^2_{yy} g(x_t, y_t-k\delta_{\epsilon} v_t) dk \big)\Big]v_t
 \nonumber \\
 & = \mathcal{S}_{[\mu,L_g]}\big[\nabla^2_{yy}g(x_t,y_t)\big]v_t.
\end{align}
Then we can obtain $\lim_{\delta_\epsilon\rightarrow 0 }\widetilde{\nabla}_v R(x_t,y_t,v_t)=\nabla_v R(x_t,y_t,v_t)$.

\section{Convergence Analysis}
In the section, we study convergence properties of our HJFBiO algorithm under some mild assumptions.
Given $x_t$ from Algorithm~\ref{alg:1}, we define a useful gradient mapping
\begin{align}  \label{eq:26}
\mathcal{G}(x_t,\nabla F(x_t),\gamma)=\frac{1}{\gamma}(x_t-x_{t+1}),
\end{align}
where $F(x)=f(x,y^*(x))$ with $y^*(x)\in \arg\min_y g(x,y)$, and $x_{t+1}$
is generated from
\begin{align}
& x_{t+1} = \mathbb{P}_{\phi(\cdot)}^\gamma(x_t ,\nabla F(x_t)) \nonumber \\
&= \arg\min_{x\in \mathbb{R}^d}\Big\{ \langle \nabla F(x_t), x\rangle
+ \frac{1}{2\gamma}\|x-x_t\|^2 + \phi(x)\Big\}. \nonumber
\end{align}
When $\phi(x)\equiv 0$, based on~(\ref{eq:26}), we have $x_{t+1}=x_t-\gamma \nabla F(x_t)$, and
then we obtain $\mathcal{G}(x_t,\nabla F(x_t),\gamma)=\nabla F(x_t)$.
\subsection{ Convergence Properties of Our Algorithm on Unimodal $g(x,y)$}
In the subsection, we study the convergence properties of our HJFBiO algorithm when $g(x,\cdot)$ satisfies the global PL condition for all $x \in \mathbb{R}^d$, i.e., it has \textbf{a unique minimizer} $y^*(x)=\arg\min_y g(x,y)$.
We first give three useful lemmas.

\begin{lemma} \label{lem:4}
Suppose the sequence $\{x_t,y_t,v_t\}_{t=1}^T$ be generated from Algorithm~\ref{alg:1}.
Under the above Assumptions, given $0<\tau \leq \frac{1}{6L_g}$, we have
\begin{align}
 & \|v_{t+1} - v^*_{t+1}\|^2 \leq ( 1-\frac{\mu\tau}{4})\|v_t -v^*_t\|^2
 - \frac{3}{4} \|v_{t+1}-v_t\|^2   \nonumber \\
 & \quad + \frac{25\tau L^2_{gyy}r^4_v\delta^2_\epsilon}{6\mu}
  + \frac{20}{3}\big(\frac{L^2_f}{\mu^2}+ \frac{L^2_{gyy}C^2_{fx}}{\mu^4}\big)\Big(\|x_{t+1}-x_t\|^2 \nonumber \\
 & \quad + \|y_{t+1}-y_t\|^2\Big),
\end{align}
where $v^*_t = v^*(x_t,y_t)=\Big(\mathcal{S}_{[\mu,L_g]}\big[\nabla^2_{yy}g(x_t,y_t)\big]\Big)^{-1}
\nabla_yf(x_t,y_t)$ for all $t\geq 1$.
\end{lemma}

\begin{lemma}\label{lem:5}
Assume the sequence $\{x_t,y_t,v_t\}_{t=1}^T$ be generated from Algorithm~\ref{alg:1}, given $0<\gamma\leq \frac{1}{2L_F}$, we have
\begin{align*}
    \Phi(x_{t+1}) &\leq \Phi(x_t) - \frac{\gamma}{2}\|\mathcal{G}(x_t,w_t,\gamma)\|^2  \nonumber \\
    & \quad + \frac{12\gamma}{\mu}\big(L^2_f+ r^2_vL^2_{gxy}\big)\big(g(x_t,y_t)-G(x_t)\big) \nonumber \\
    & \quad + 6\gamma C^2_{gxy}\|v_t-v^*_t\|^2 + 2\gamma L^2_{gxy}\delta^2_{\epsilon}r^4_v,
\end{align*}
where $\Phi(x) = F(x)+\phi(x)$ and $\mathcal{G}(x_t,w_t,\gamma) = \frac{1}{\gamma}(x_t-x_{t+1})$.
\end{lemma}

\begin{lemma} \label{lem:6}
Suppose the sequence $\{x_t,y_t,v_t\}_{t=1}^T$ be generated from Algorithm~\ref{alg:1}.
Under the above Assumptions~\ref{ass:1}-\ref{ass:3}, given $\gamma\leq \min\Big\{\frac{\lambda\mu }{16L_G }, \frac{\mu }{16L^2_g } \Big\}$ and $0<\lambda \leq \frac{1}{2L_g }$, we have
\begin{align}
& g(x_{t+1},y_{t+1}) - G(x_{t+1}) \nonumber \\
& \leq (1-\frac{\lambda\mu}{2 }) \big(g(x_t,y_t) -G(x_t)\big) + \frac{1}{8\gamma}\|x_{t+1}-x_t\|^2 \nonumber \\
& \quad   -\frac{1}{4\lambda }\|y_{t+1}-y_t\|^2 + \lambda\|\nabla_y g(x_t,y_t)-u_t\|^2,
\end{align}
where $G(x_t)=g(x_t,y^*(x_t))$ with $y^*(x_t) \in \arg\min_{y}g(x_t,y)$ for all $t\geq 1$.
\end{lemma}

Based on the above useful lemmas, we give the convergence properties of our HJFBiO method in the following.
\begin{theorem}  \label{th:1}
 Assume the sequence $\{x_t,y_t,v_t\}_{t=1}^T$ be generated from our Algorithm \ref{alg:1}. Under the above Assumptions~\ref{ass:1}-\ref{ass:5}, let $0<\gamma\leq \min\Big(\frac{1}{2L_F},\frac{\lambda\mu }{16L_G}, \frac{\mu }{16L^2_g},\frac{3}{160\breve{L}^2}, \frac{\mu\tau}{30C^2_{gxy}},\frac{\mu^2\lambda}{30(L^2_f+ r^2_vL^2_{gxy})}\Big)$, $0< \lambda\leq \min\big(\frac{1}{2L_g},\frac{3}{80\breve{L}^2}\big)$ and $0<\tau \leq \frac{1}{6L_g}$,
 we have
 \begin{align}
 & \frac{1}{T}\sum_{t=1}^T\|\mathcal{G}(x_t,\nabla F(x_t),\gamma)\|^2 \nonumber \\
 & \leq \frac{8(\Phi(x_1) +  g(x_1,y_1)-G(x_1)+ \|v_1-v_1^*\|^2 - \Phi^*)}{T\gamma} \nonumber \\
 & \quad + 20 L^2_{gxy}\delta^2_{\epsilon}r^4_v + \frac{100\tau L^2_{gyy}r^4_v\delta^2_\epsilon}{3\gamma\mu},
\end{align}
where $\breve{L}^2=\frac{L^2_f}{\mu^2}+ \frac{L^2_{gyy}C^2_{fx}}{\mu^4}$.
\end{theorem}

\begin{remark}
Without loss of generality, let $\tau = O(1)$, $\lambda=O(1)$ and $\gamma=\min\big(\frac{1}{2L_F},\frac{\lambda\mu }{16L_G}, \frac{\mu }{16L^2_g},\frac{3}{160\breve{L}^2},\\ \frac{\mu\tau}{30C^2_{gxy}},\frac{\mu^2\lambda}{30(L^2_f+ r^2_vL^2_{gxy})} \big)=O(1)$.
Furthermore, let  $\delta_{\epsilon}=O\big(\frac{1}{\sqrt{T}\max(L^2_{gxy},L^2_{gyy}/\mu)r^2_v}\big)$, we have
\begin{align*}
\frac{1}{T}\sum_{t=1}^T\|\mathcal{G}(x_t,\nabla F(x_t),\gamma)\|^2  \leq O(\frac{1}{T}).
\end{align*}
Let $O(\frac{1}{T}) = \epsilon$, we obtain $T=(\epsilon^{-1})$. Since requiring seven first-order gradients at each iteration, our HJFBiO algorithm can obtain an optimal gradient (or iteration) complexity of $7\cdot T=O(\epsilon^{-1})$ in finding an $\epsilon$-stationary solution of Problem~(\ref{eq:1}), which matches the lower bound established by the first-order method for finding an $\epsilon$-stationary point of nonconvex smooth optimization~\citep{carmon2020lower}.

When $\phi(x)\equiv 0$, based on~(\ref{eq:26}), we have $\mathcal{G}(x_t,\nabla F(x_t),\gamma)=\nabla F(x_t)$.
Thus our HJFBiO algorithm still obtains an optimal gradient complexity of $7\cdot T=O(\epsilon^{-1})$ in finding an $\epsilon$-stationary solution of Problem~(\ref{eq:1}) (i.e., $\min_{1\leq t\leq T}\|\nabla F(x_t)\|^2\leq \epsilon$).
\end{remark}
\subsection{Convergence Properties of Our Algorithm on multimodal $g(x,y)$} \label{sec4.2}
In this subsection, we study the convergence properties of our HJFBiO method when $g(x,\cdot)$ satisfies the local PL condition for all $x$, i.e, it has \textbf{multi local minimizers } $y^\diamond(x,y)\in \arg\min_y g(x,y)$. As in~\cite{liu2022bome}, we define the attraction points.
\begin{definition}
Given any $(x,y)$, if sequence $\{y_t\}_{t=0}^{\infty}$ generated by gradient descent $y_t=y_{t-1}-\lambda \nabla_y g(x,y_{t-1})$ starting from $y_0=y$ converges to $y^\diamond(x,y)$, we say that $y^\diamond(x,y)$ is
the attraction point of $(x,y)$ with step size $\lambda>0$.
\end{definition}
An attraction basin be formed by the same attraction point in set of $(x,y)$. In the following analysis,
we assume the PL condition within the individual attraction basins. Let $F^\diamond(x)=f(x,y^\diamond(x,y))$.

\begin{assumption} \label{ass:1g}
(\textbf{Local PL Condition in Attraction Basins})
Assume that for any $(x,y)$, $y^\diamond(x,y)$ exists.
 $g(x,\cdot)$ satisfies the local PL condition in attraction basins,
i.e., for any $(x,y)$, there exists a constant $\mu>0$ such that
\begin{align}
 \|\nabla_y g(x,y)\|^2 \geq 2\mu \big(g(x,y)- G^\diamond(x)\big),
\end{align}
where $G^\diamond(x)=g(x,y^\diamond(x,y))$.
\end{assumption}
\begin{assumption} \label{ass:2g}
 The function $g\big(x,y^\diamond(x,y)\big)$ satisfies
\begin{align}
 \varrho\big(\nabla^2_{yy}g\big(x,y^\diamond(x,y)\big)\big) \in [\mu, L_g],
\end{align}
where $y^\diamond(x,y)$ is the attraction point of $(x,y)$, and $\varrho(\cdot)$ denotes the eigenvalue (or singular-value) function and $L_g\geq\mu>0$.
\end{assumption}
\begin{assumption} \label{ass:5g}
 The function $\Phi^\diamond(x)=F^\diamond(x)+\phi(x)$ is bounded below in $\mathbb{R}^d$, \emph{i.e.,} $\Phi^\diamond = \inf_{x\in \mathbb{R}^d} \Phi^\diamond(x) > -\infty$.
\end{assumption}

\begin{theorem}  \label{th:1g}
 Assume the sequence $\{x_t,y_t,v_t\}_{t=1}^T$ be generated from our Algorithm \ref{alg:1}. Under the above Assumptions~\ref{ass:1g},~\ref{ass:2g},~\ref{ass:3},~\ref{ass:4},~\ref{ass:5g}, let $0<\gamma\leq \min\Big(\frac{1}{2L_F},\frac{\lambda\mu }{16L_G}, \frac{\mu }{16L^2_g},\frac{3}{160\breve{L}^2}, \frac{\mu\tau}{30C^2_{gxy}},\frac{\mu^2\lambda}{30(L^2_f+ r^2_vL^2_{gxy})}\Big)$, $0< \lambda\leq \min\big(\frac{1}{2L_g},\frac{3}{80\breve{L}^2}\big)$ and $0<\tau \leq \frac{1}{6L_g}$,
 we have
 \begin{align}
 & \frac{1}{T}\sum_{t=1}^T\|\mathcal{G}(x_t,\nabla F^\diamond(x_t),\gamma)\|^2 \nonumber \\
 & \leq \frac{8(\Phi^\diamond(x_1) +  g(x_1,y_1)-G(x_1)+ \|v_1-v_1^*\|^2 - \Phi^\diamond)}{T\gamma} \nonumber \\
 & \quad + 20 L^2_{gxy}\delta^2_{\epsilon}r^4_v + \frac{100\tau L^2_{gyy}r^4_v\delta^2_\epsilon}{3\gamma\mu},
\end{align}
where $\Phi^\diamond(x)=F^\diamond(x)+\phi(x)$ and $\breve{L}^2=\frac{L^2_f}{\mu^2}+ \frac{L^2_{gyy}C^2_{fx}}{\mu^4}$.
\end{theorem}
\begin{remark}
The proof of Theorem~\ref{th:1g} \textbf{can follow} the proof of Theorem~\ref{th:1}.
Let further $\delta_{\epsilon}=O\big(\frac{1}{\sqrt{T}\max(L^2_{gxy},L^2_{gyy}/\mu)r^2_v}\big)$,
our HJFBiO algorithm can also obtain an optimal gradient (or iteration) complexity of $7\cdot T=O(\epsilon^{-1})$ in finding an $\epsilon$-stationary solution of Problem~(\ref{eq:1}) \textbf{under local PL condition}.
\end{remark}

\section{ Experiments}
In the section, we conduct bilevel PL game and
hyper-representation learning tasks
to demonstrate efficiency of our method. In the experiments, we compare our method
with the existing methods given in Table~\ref{tab:1}. Meanwhile, we add a baseline: BVFSM~\citep{liu2021value}.
For fair comparison, since only the AdaPAG~\citep{huang2023adaptive} uses adaptive learning rate, we
exclude it in the comparisons.
\begin{figure}[ht]
\vspace*{-8pt}
\centering
 \subfloat{\includegraphics[width=0.24\textwidth]{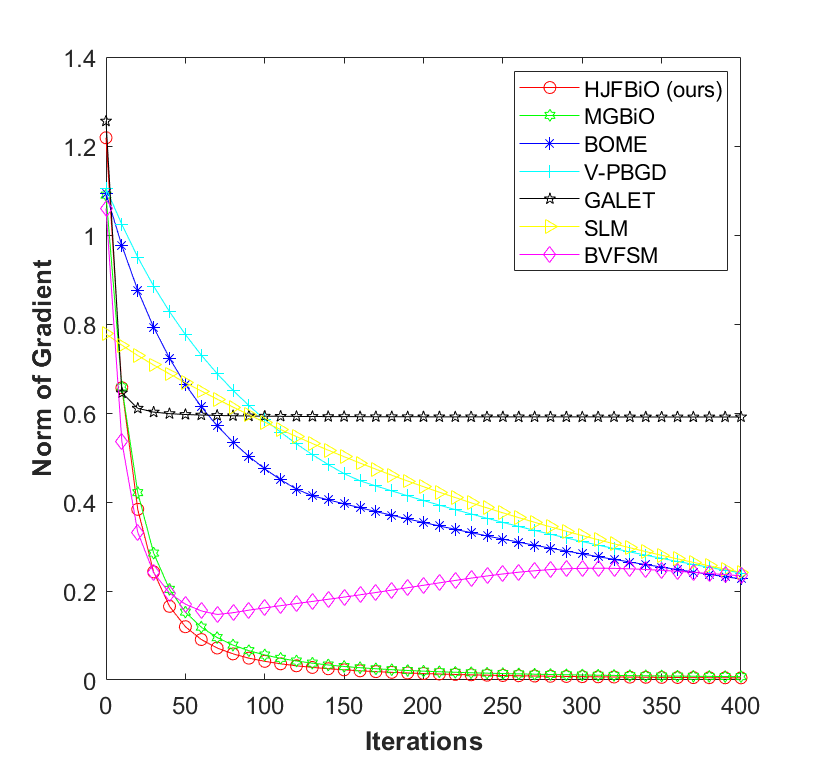}}
  \hfill
 \subfloat{\includegraphics[width=0.24\textwidth]{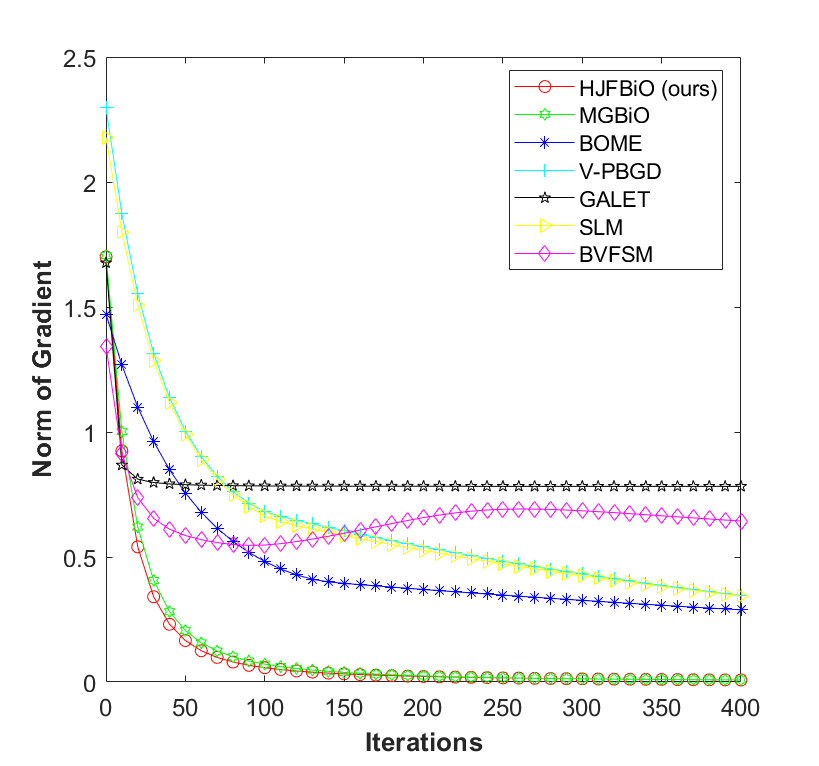}}
  \hfill
\caption{ PL Game: norm of gradient vs number of iteration under $d=100$ (Left) and $d=200$ (Right).}
\label{fig:1}
\end{figure}

\subsection{ Bilevel Polyak-{\L}ojasiewicz Game}
In this subsection, as in \citep{huang2023momentum},
we apply the bilevel Polyak-{\L}ojasiewicz game task to verify efficiency of our algorithm,
defined as:
\begin{align}
 \min_{x\in \mathbb{R}^d} &\ \frac{1}{2}x^TPx + x^TR^1y,   \\
 \mbox{s.t.} & \ \min_{y\in \mathbb{R}^d} \frac{1}{2}y^TQy + x^TR^2y, \nonumber
\end{align}
where $P=\frac{1}{n}\sum_{i=1}^n p_i(p_i)^T$, $Q=\frac{1}{n}\sum_{i=1}^n q_i(q_i)^T$,
$R^1= \frac{1}{n}\sum_{i=1}^n 0.01r^1_i(r^1_i)^T$ and $R^2= \frac{1}{n}\sum_{i=1}^n 0.01r^2_i(r^2_i)^T$.
Specifically, samples $\{p_i\}_{i=1}^n$,
$\{q_i\}_{i=1}^n$, $\{r^1_i\}_{i=1}^n$ and $\{r^2_i\}_{i=1}^n$ are independently drawn from normal distributions $\mathbb{N}(0,\Sigma_{P})$, $\mathbb{N}(0,\Sigma_{Q})$, $\mathbb{N}(0,\Sigma_{R^1})$ and $\mathbb{N}(0,\Sigma_{R^2})$, respectively.
Here we set $\Sigma_{P}=U^1D^1(U^1)^T$, where $U^1\in \mathbb{R}^{d\times l} \ (l<d)$
is column orthogonal, and $D^1\in \mathbb{R}^{l\times l}$ is diagonal and its diagonal elements are distributed uniformly in the interval $[\mu,L]$ with $0<\mu<L$. Let $\Sigma_{Q}=U^2D^2(U^2)^T$,
where $U^2\in \mathbb{R}^{d\times l}$
is column orthogonal, and $D^2\in \mathbb{R}^{l\times l}$ is diagonal and its diagonal elements are distributed uniformly in the interval $[\mu,L]$ with $0<\mu<L$.
We also set $\Sigma_{R^1}=0.001V^1(V^1)^T$ and $\Sigma_{R^2}=0.001V^2(V^2)^T$, where each element of $V^1, V^2\in \mathbb{R}^{d\times d}$ is independently sampled from normal distribution $\mathbb{N}(0,1)$. Since
the covariance matrices $\Sigma_P$ and $\Sigma_Q$ are rank-deficient,
it is ensured that both $P$ and $Q$ are singular.

In the experiment, we set $l=50$, $n=50\cdot d$, $d=100$ and $d=200$. For fair comparison, we set a basic learning rate as $0.01$ for all algorithms. In our HJFBiO method, we set $\delta_\epsilon = 10^{-5}$.  Figure~\ref{fig:1} provides the results on \emph{norm of gradient} vs \emph{iteration}, where the iteration denotes iteration at outer loop in all algorithms.  Here these results verify that our HJFBiO algorithm outperforms all comparisons. In particular, our HJFBiO method has a lower computation at each iteration than the MGBiO method.

\begin{figure}[ht]
\vspace*{-8pt}
\centering
 \subfloat{\includegraphics[width=0.24\textwidth]{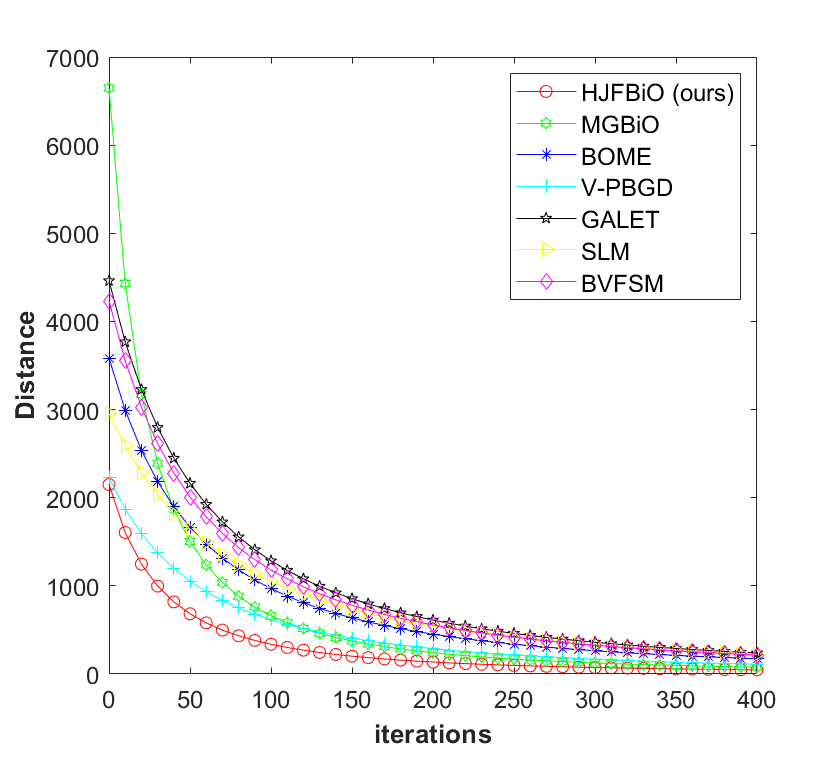}}
  \hfill
 \subfloat{\includegraphics[width=0.24\textwidth]{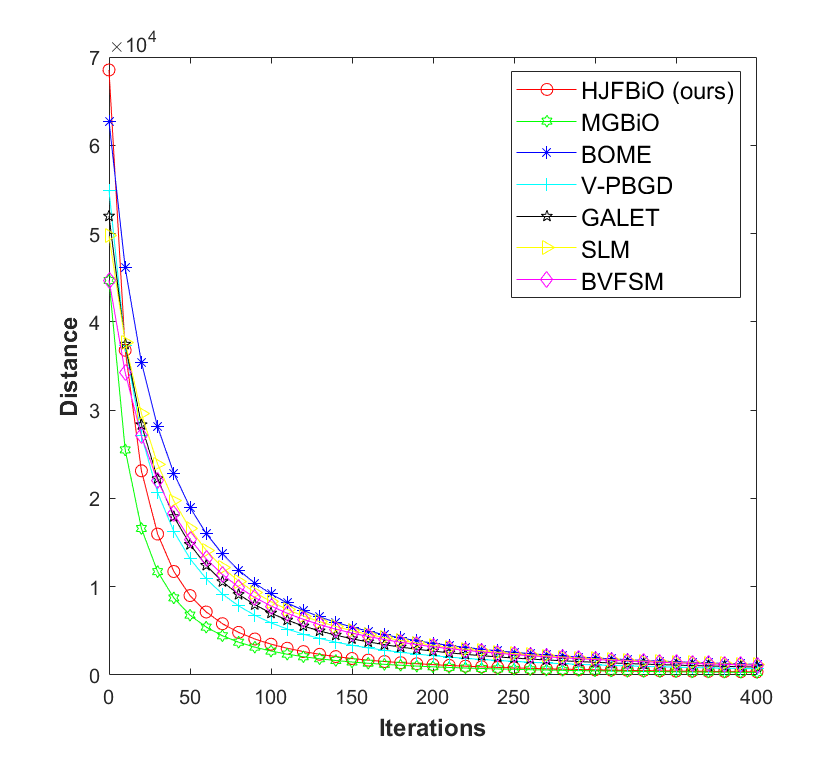}}
  \hfill
\caption{Distances of the algorithms under the case of $d=100$ (Left) and $d=200$ (Right).}
\label{fig:2}
\end{figure}

\subsection{ Hyper-Representation Learning}
In this subsection, as in \citep{huang2023momentum}, we consider the hyper-representation learning on matrix sensing task to verify efficiency of our method.
Specifically, given $n$ sensing matrices $\{C_i\in \mathbb{R}^{d\times d}\}_{i=1}^n$ and $n$ observations $o_i=\langle C_i, H^*\rangle=\mbox{trace}(C_i^TH^*)$, where $H^*=U^*(U^*)^T$ is a low-rank symmetric matrix with
$U^*\in \mathbb{R}^{d\times r}$, the goal of this task is to find an optimal matrix $U^*$, which can be
defined as the following problem:
\begin{align}
 \min_{U\in \mathbb{R}^{d\times r}} \frac{1}{n}\sum_{i=1}^n \ell_i(U)\equiv \frac{1}{2}\big(\langle C_i, UU^T\rangle-o_i\big)^2.
\end{align}
Next, we consider the hyper-representation learning in matrix sensing task,
 which be rewritten the following
bilevel optimization problem:
\begin{align}
& \min_{x\in \mathbb{R}^{d\times r-1}}\frac{1}{|D_v|}\sum_{i\in D_v} \ell_i(x,y^*(x)),  \\
 & \mbox{s.t.} \  y^*(x)\in \arg\min_{y\in \mathbb{R}^{d}}\frac{1}{|D_t|}\sum_{i\in D_t} \ell_i(x,y), \nonumber
\end{align}
where $U=[y;x]\in \mathbb{R}^{d\times r}$ is a concatenation of $x$ and $y$. Here we define variable $x$ to be the first $r-1$ columns of $U$ and variable $y$ to be the last column. Meanwhile, $D_t$ denotes the training dataset, and $D_v$ denotes the validation dataset.
The ground truth low-rank matrix $H^*$ is generated by $H^*=U^*(U^*)^T$, where each entry of $U^*$ is drawn from
normal distribution $\mathbb{N}(0,1/d)$ independently. We randomly generate $n=30\cdot d$ samples of sensing
matrices $\{C_i\}_{i=1}^n$ from standard normal distribution, and then compute the corresponding no-noise
labels $o_i = \langle C_i, H^*\rangle$. We split all samples into two dataset: a train dataset
$D_t$ with 40\% data and a validation dataset $D_v$ with 60\% data.

\begin{figure}[ht]
\vspace*{-8pt}
\centering
 \subfloat{\includegraphics[width=0.24\textwidth]{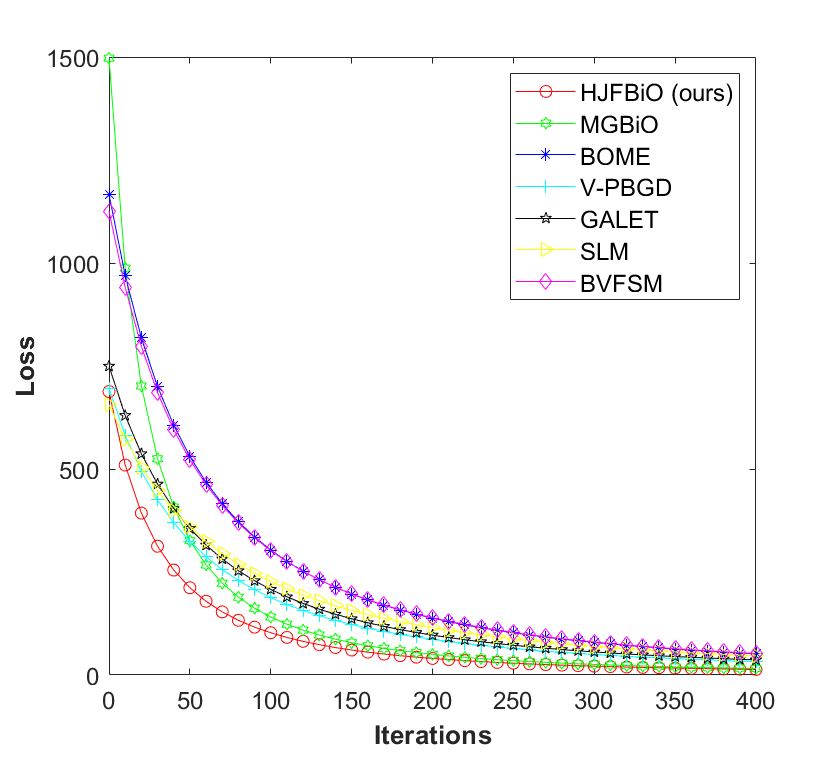}}
  \hfill
 \subfloat{\includegraphics[width=0.24\textwidth]{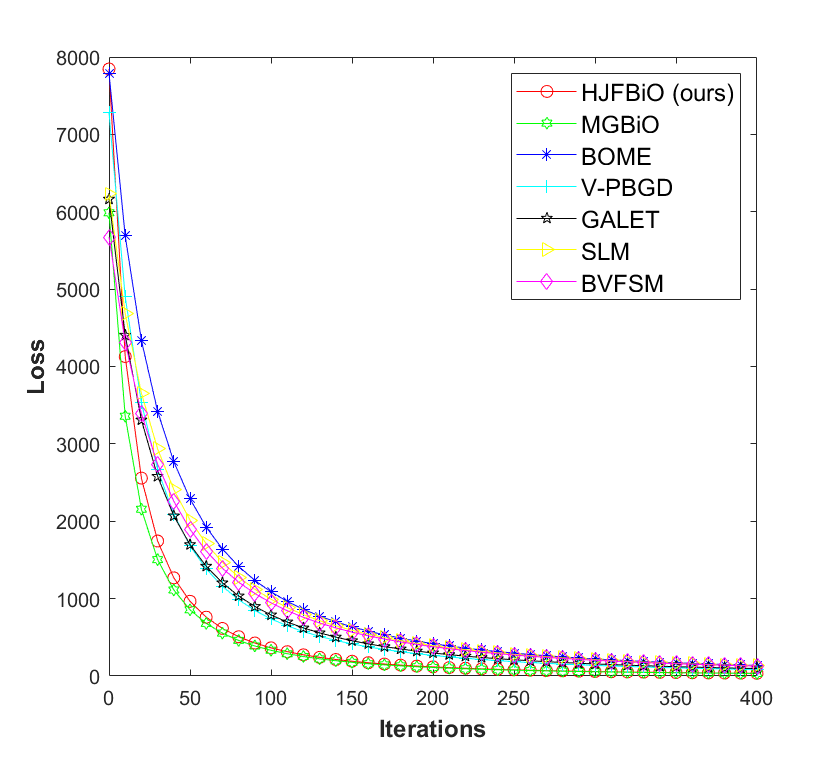}}
  \hfill
\caption{Losses of the algorithms under the case of $d=100$ (Left) and $d=200$ (Right).}
\label{fig:3}
\end{figure}

In the experiment, for fair comparison, we set the basic learning rate as
$0.01$ for all algorithms. In our HJFBiO method, we set $\delta_\epsilon = 10^{-5}$.
Let $\ell(U) = \frac{1}{2n}\sum_{i=1}^n\big(\langle C_i, UU^T\rangle-o_i\big)^2$ denote the loss, and
$\|UU^T-H^*\|^2_F/\|H^*\|^2_F$ denotes the distance.
Figures~\ref{fig:2} and~\ref{fig:3} show that our HJFBiO method outperforms all the comparisons
on the \textbf{distance} vs iteration, where the iteration denotes iteration at outer loop in all algorithms.
While our HJFBiO method is comparable with the MGBiO method on the \textbf{loss} vs iteration.
In particular, our HJFBiO method has a lower computation at each iteration than the MGBiO method.

\section{Conclusions }
In the paper, we proposed an efficient Hessian/Jacobian-free bilevel method to solve the nonconvex-PL bilevel problems based on the finite-difference estimator and a new projection operator.
Moreover, under some mild assumptions, we proved that our HJFBiO method obtains the best known convergence rate $O(\frac{1}{T})$, and gets an optimal gradient (or iteration) complexity of $O(\epsilon^{-1})$ in finding $\epsilon$-stationary solution under global and local PL condition, respectively.

%% Acknowledgements should only appear in the accepted version.
\section*{Acknowledgements}
We thank the anonymous reviewers for their helpful comments. This paper was partially supported by NSFC under Grant No.
62376125. It was also partially supported by the Fundamental Research Funds for the Central Universities NO.NJ2023032.

\section*{Impact Statement}
This paper presents work whose goal is to advance the field of Machine Learning. There are many potential societal consequences of our work, none which we feel must be specifically highlighted here.

\nocite{langley00}

\bibliography{HFBIO}
\bibliographystyle{icml2024}

%%%%%%%%%%%%%%%%%%%%%%%%%%%%%%%%%%%%%%%%%%%%%%%%%%%%%%%%%%%%%%%%%%%%%%%%%%%%%%%
%%%%%%%%%%%%%%%%%%%%%%%%%%%%%%%%%%%%%%%%%%%%%%%%%%%%%%%%%%%%%%%%%%%%%%%%%%%%%%%
% APPENDIX
%%%%%%%%%%%%%%%%%%%%%%%%%%%%%%%%%%%%%%%%%%%%%%%%%%%%%%%%%%%%%%%%%%%%%%%%%%%%%%%
%%%%%%%%%%%%%%%%%%%%%%%%%%%%%%%%%%%%%%%%%%%%%%%%%%%%%%%%%%%%%%%%%%%%%%%%%%%%%%%
\newpage
\appendix
\onecolumn

\section{Detailed Convergence Analysis}
In this section, we provide the detailed convergence analysis of our algorithms.
We first review some useful lemmas.

\begin{lemma} \label{lem:A1}
(\cite{nesterov2018lectures})
Assume that $f(x)$ is a differentiable convex function and $\mathcal{X}$ is a convex set.
 $x^* \in \mathcal{X}$ is the solution of the
constrained problem $\min_{x\in \mathcal{X}}f(x)$, if
\begin{align}
 \langle \nabla f(x^*), x-x^*\rangle \geq 0, \ \forall x\in \mathcal{X}.
\end{align}
\end{lemma}

\begin{lemma} \label{lem:A2}
(\cite{karimi2016linear})
 The function $f(x): \mathbb{R}^d\rightarrow \mathbb{R}$ is $L$-smooth and satisfies PL condition with constant $\mu$, then it also satisfies error bound (EB) condition with $\mu$, i.e., for all $x \in \mathbb{R}^d$
\begin{align}
 \|\nabla f(x)\| \geq \mu\|x^*-x\|,
\end{align}
where $x^* \in \arg\min_{x} f(x)$. It also satisfies quadratic growth (QG) condition with $\mu$, i.e.,
\begin{align}
 f(x)-\min_x f(x) \geq \frac{\mu}{2}\|x^*-x\|^2.
\end{align}
\end{lemma}

\subsection{Convergence Analysis of HJFBiO Algorithm for Bilevel Optimization with Regularization}
In this subsection, we detail the convergence analysis of our HJFBiO algorithm
for bilevel optimization.
We give some useful lemmas.

\begin{lemma} \label{lem:B1}
(Restatement of Lemma~\ref{lem:4})
Suppose the sequence $\{x_t,y_t,v_t\}_{t=1}^T$ be generated from Algorithm~\ref{alg:1}.
Under the above Assumptions~\ref{ass:1}-\ref{ass:4}, given $0<\tau \leq \frac{1}{6L_g}$, we have
\begin{align}
 \|v_{t+1} - v^*_{t+1}\|^2 & \leq ( 1-\frac{\mu\tau}{4})\|v_t -v^*_t\|^2
 - \frac{3}{4} \|v_{t+1}-v_t\|^2   \nonumber \\
 & \quad + \frac{25\tau L^2_{gyy}r^4_v\delta^2_\epsilon}{6\mu}
  + \frac{20}{3}\big(\frac{L^2_f}{\mu^2}+ \frac{L^2_{gyy}C^2_{fx}}{\mu^4}\big)\big(\|x_{t+1}-x_t\|^2 + \|y_{t+1}-y_t\|^2\big),
\end{align}
where $v^*_t = v^*(x_t,y_t)=\Big(\mathcal{S}_{[\mu,L_g]}\big[\nabla^2_{yy}g(x_t,y_t)\big]\Big)^{-1}
\nabla_yf(x_t,y_t)$ for all $t\geq 1$.
\end{lemma}

\begin{proof}
 Since the function $ R(x,y,v) = \frac{1}{2}v^T\mathcal{S}_{[\mu,L_g]}\big[\nabla^2_{yy}g(x,y)\big]v - v^T\nabla_yf(x,y)$ is $\mu$-strongly convex on variable $v$,
we have
\begin{align} \label{eq:E1}
 R(x_t,y_t,v) & \geq R(x_t,y_t,v_t) + \langle\nabla_v R(x_t,y_t,v_t), v-v_t\rangle + \frac{\mu}{2}\|v-v_t\|^2 \nonumber \\
 & = R(x_t,y_t,v_t) + \langle h_t, v-v_{t+1}\rangle + \langle\nabla_v R(x_t,y_t,v_t)-h_t, v-v_{t+1}\rangle \nonumber \\
 & \quad +\langle\nabla_v R(x_t,y_t,v_t), v_{t+1}-v_t\rangle + \frac{\mu}{2}\|v-v_t\|^2.
\end{align}
Since the function $R(x,y,v)$ is $L_g$-smooth on variable $v$, we have
\begin{align} \label{eq:E2}
 R(x_t,y_t,v_{t+1}) \leq R(x_t,y_t,v_t) + \langle\nabla_v R(x_t,y_t,v_t), v_{t+1}-v_t\rangle + \frac{L_g}{2}\|v_{t+1}-v_t\|^2 .
\end{align}
By combining the about inequalities \eqref{eq:E1} with \eqref{eq:E2}, we have
\begin{align} \label{eq:E3}
 R(x_t,y_t,v) & \geq R(x_t,y_t,v_{t+1}) + \langle h_t, v-v_{t+1}\rangle + \langle\nabla_v R(x_t,y_t,v_t)-h_t, v-v_{t+1}\rangle \nonumber \\
 & \quad + \frac{\mu}{2}\|v-v_t\|^2 - \frac{L_g}{2}\|v_{t+1}-v_t\|^2.
\end{align}

According to the line 5 of Algorithm~\ref{alg:1}, we have $v_{t+1}=\mathcal{P}_{r_v}\big(v_t - \tau h_t\big) = \arg\min_{v\in \Lambda_v } \Big\{\langle h_t,v-v_t\rangle + \frac{1}{2\tau}\|v-v_t\|^2\Big\}$ with $\Lambda_v=\{v\in \mathbb{R}^p \ | \|v\|\leq r_v\}$. According to the above Lemma~\ref{lem:A1}, given $v\in \Lambda_v$, we have
\begin{align} \label{eq:E4}
  \langle h_t + \frac{1}{\tau}(v_{t+1}-v_t), v-v_{t+1}\rangle \geq 0.
 \end{align}

By plugging the inequalities \eqref{eq:E4} into \eqref{eq:E3}, we have
\begin{align}
 R(x_t,y_t,v) & \geq R(x_t,y_t,v_{t+1}) + \frac{1}{\tau}\langle v_{t+1}- v_t, v_t- v\rangle + \frac{1}{\tau}\|v_{t+1}- v_t\|^2 \nonumber \\
 & \quad +  \langle\nabla_v R(x_t,y_t,v_t)-h_t, v-v_{t+1}\rangle + \frac{\mu}{2}\|v-v_t\|^2 - \frac{L_g}{2}\|v_{t+1}-v_t\|^2.
\end{align}
Let $v=v^*_t = v^*(x_t,y_t)=\Big(\mathcal{S}_{[\mu,L_g]}\big[\nabla^2_{yy}g(x_t,y_t)\big]\Big)^{-1}
\nabla_yf(x_t,y_t)$, then we have
\begin{align}
 R(x_t,y_t,v^*_t) & \geq R(x_t,y_t,v_{t+1}) + \frac{1}{\tau}\langle v_{t+1}- v_t, v_t - v^*_t\rangle
  + (\frac{1}{\tau}-\frac{L_g}{2})\|v_{t+1}- v_t\|^2 \nonumber \\
 & \quad + \langle\nabla_v R(x_t,y_t,v_t)-h_t, v^*_t-v_{t+1}\rangle + \frac{\mu}{2}\|v^*_t-v_t\|^2.
\end{align}
Due to the strongly-convexity of $R(x,y,\cdot)$ and $v^*_t =\arg\min_{v\in \Lambda_v}R(x_t,y_t,v)$ with $r_v=\frac{C_{fy}}{\mu}$, we have $R(x_t,y_t,v^*_t) \leq R(x_t,y_t,v_{t+1})$.
Thus, we obtain
\begin{align} \label{eq:E5}
 0 & \geq  \frac{1}{\tau}\langle v_{t+1}- v_t, v_t - v^*_t\rangle
  + (\frac{1}{\tau}-\frac{L_g}{2})\|v_{t+1}- v_t\|^2 \nonumber \\
 & \quad + \langle\nabla_v R(x_t,y_t,v_t)-h_t, v^*_t-v_{t+1}\rangle + \frac{\mu}{2}\|v^*_t-v_t\|^2.
\end{align}

Consider the term $\langle v_{t+1}- v_t, v_t - v^*_t\rangle$, we have
\begin{align} \label{eq:E6}
 \langle v_{t+1}- v_t, v_t - v^*_t\rangle = \frac{1}{2}\|v_{t+1}-v^*_t\|^2 - \frac{1}{2}\|v_t -v^*_t\|^2 - \frac{1}{2}\|v_{t+1}-v_t\|^2.
\end{align}
Consider the upper bound of the term $\langle\nabla_v R(x_t,y_t,v_t)-h_t, v^*_t-v_{t+1}\rangle$, we have
\begin{align} \label{eq:E7}
 &\langle\nabla_v R(x_t,y_t,v_t)-h_t, v^*_t-v_{t+1}\rangle \nonumber \\
 & = \langle\nabla_v R(x_t,y_t,v_t)-h_t, v^*_t-v_t\rangle + \langle\nabla_v R(x_t,y_t,v_t)-h_t, v_t-v_{t+1}\rangle \nonumber \\
 & \geq -\frac{1}{\mu} \|\nabla_v R(x_t,y_t,v_t)-h_t\|^2 - \frac{\mu}{4}\|v^*_t-v_t\|^2 - \frac{1}{\mu} \|\nabla_v R(x_t,y_t,v_t)-h_t\|^2 - \frac{\mu}{4}\|v_t-v_{t+1}\|^2 \nonumber \\
 & = -\frac{2}{\mu} \|\nabla_v R(x_t,y_t,v_t)-h_t\|^2 - \frac{\mu}{4}\|v^*_t-v_t\|^2 - \frac{\mu}{4}\|v_t-v_{t+1}\|^2 \nonumber\\
 & \geq -\frac{2L^2_{gyy}r^4_v\delta^2_\epsilon}{\mu}- \frac{\mu}{4}\|v^*_t-v_t\|^2 - \frac{\mu}{4}\|v_t-v_{t+1}\|^2,
\end{align}
where the last inequality holds by the following inequality:
\begin{align}
 & \|\nabla_v R(x_t,y_t,v_t)-h_t\|^2 \nonumber \\
 & = \|\nabla_v R(x_t,y_t,v_t)-\widetilde{\nabla}_v R(x_t,y_t,v_t)\|^2 \nonumber \\
 & = \big\|\mathcal{S}_{[\mu,L_g]}\big[\nabla^2_{yy}g(x_t,y_t)\big]v_t - \nabla_yf(x_t,y_t)-\mathcal{M}_{r_h}\big(\widetilde{H}(x_t, y_t, v_t, \delta_{\epsilon})\big) + \nabla_yf(x_t,y_t)\big\|^2 \nonumber \\
 & = \big\|\mathcal{S}_{[\mu,L_g]}\big[\nabla^2_{yy}g(x_t,y_t)\big]v_t-\mathcal{M}_{r_h}\big(\widetilde{H}(x_t, y_t, v_t, \delta_{\epsilon})\big)\big\|^2   \nonumber \\
 & = \big\|\mathcal{S}_{[\mu,L_g]}\big[\nabla^2_{yy}g(x_t,y_t)\big]v_t-\mathcal{M}_{r_h}\Big(\frac{1}{2\delta_\epsilon}\big(\nabla_y g(x_t, y_t+\delta_{\epsilon} v_t) - \nabla_y g(x_t, y_t-\delta_{\epsilon} v_t)\big)\Big)\big\|^2   \nonumber \\
  & = \big\|\mathcal{S}_{[\mu,L_g]}\big[\nabla^2_{yy}g(x_t,y_t)\big]v_t- \mathcal{M}_{r_h}\Big(\frac{1}{2\delta_\epsilon}\big(\nabla_y g(x_t, y_t+\delta_{\epsilon} v_t) -\nabla_y g(x_t, y_t) + \nabla_y g(x_t, y_t)- \nabla_y g(x_t, y_t-\delta_{\epsilon} v_t)\big)\Big)\big\|^2   \nonumber \\
 & = \big\|\mathcal{S}_{[\mu,L_g]}\big[\nabla^2_{yy}g(x_t,y_t)\big]v_t- \mathcal{M}_{r_h}\Big(
\frac{1}{2\delta_{\epsilon}} \int_{k=0}^1\nabla^2_{yy} g(x_t, y_t+k\delta_{\epsilon} v_t) \delta_{\epsilon}v_t dk + \frac{1}{2\delta_{\epsilon}}\int_{k=0}^1\nabla^2_{yy} g(x_t, y_t-k\delta_{\epsilon} v_t) \delta_{\epsilon}v_t dk\Big)\big\|^2   \nonumber \\
 & = \big\|\mathcal{S}_{[\mu,L_g]}\big[\nabla^2_{yy}g(x_t,y_t)\big]v_t-  \mathcal{M}_{r_h}\Big(
 \big(\frac{1}{2}\int_{k=0}^1\nabla^2_{yy} g(x_t, y_t+k\delta_{\epsilon} v_t) dk + \frac{1}{2}\int_{k=0}^1\nabla^2_{yy} g(x_t, y_t-k\delta_{\epsilon} v_t) dk \big)v_t\Big) \big\|^2   \nonumber \\
 & \mathop{=}^{(i)} \big\|\mathcal{S}_{[\mu,L_g]}\big[\nabla^2_{yy}g(x_t,y_t)\big]v_t-  \mathcal{S}_{[\mu,L_g]}\Big[\frac{1}{2}\int_{k=0}^1\nabla^2_{yy} g(x_t, y_t+k\delta_{\epsilon} v_t) dk + \frac{1}{2}\int_{k=0}^1\nabla^2_{yy} g(x_t, y_t-k\delta_{\epsilon} v_t) dk \big)\Big]v_t \big\|^2   \nonumber \\
 & \mathop{\leq}^{(ii)}  r^2_v\big\|\nabla^2_{yy}g(x_t,y_t)-\frac{1}{2}\int_{k=0}^1\nabla^2_{yy} g(x_t, y_t+k\delta_{\epsilon} v_t) dk - \frac{1}{2}\int_{k=0}^1\nabla^2_{yy} g(x_t, y_t-k\delta_{\epsilon} v_t) dk \big)\big\|^2 \nonumber \\
 & \leq \frac{r^2_v}{2}\Big(\int_{k=0}^1\big\|\nabla^2_{yy}g(x_t,y_t)-\nabla^2_{yy} g(x_t, y_t+k\delta_{\epsilon} v_t)\big\| dk \Big)^2 + \frac{r^2_v}{2}\Big(\int_{k=0}^1\big\|\nabla^2_{yy}g(x_t,y_t)-\nabla^2_{yy} g(x_t, y_t-k\delta_{\epsilon} v_t)\big\| dk \Big)^2 \nonumber \\
 & \leq L^2_{gyy}r^4_v\delta^2_\epsilon,
\end{align}
where the above equality (i) holds by $v_t = \mathcal{P}_{r_v}(v_t)$ and the above definition~\ref{def:1}, and
the above inequality (ii) holds by $\|v_t\| \leq r_v$.

By plugging the inequalities \eqref{eq:E6} and \eqref{eq:E7} into \eqref{eq:E5},
we obtain
\begin{align}
  \frac{1}{2\tau}\|v_{t+1}-v^*_t\|^2
 & \leq (\frac{1}{2\tau}-\frac{\mu}{4})\|v_t -v^*_t\|^2 + (  \frac{\mu}{4} + \frac{L_g}{2} -\frac{1}{2\tau} ) \|v_{t+1}-v_t\|^2 + \frac{2L^2_{gyy}r^4_v\delta^2_\epsilon}{\mu} \nonumber \\
 & \leq (\frac{1}{2\tau}-\frac{\mu}{4})\|v_t -v^*_t\|^2 + ( \frac{3L_g}{4} -\frac{1}{2\tau} ) \|v_{t+1}-v_t\|^2 + \frac{2L^2_{gyy}r^4_v\delta^2_\epsilon}{\mu} \nonumber \\
 & = (\frac{1}{2\tau}-\frac{\mu}{4})\|v_t -v^*_t\|^2 - ( \frac{3}{8\tau} + \frac{1}{8\tau}- \frac{3L_g}{4}) \|v_{t+1}-v_t\|^2 + \frac{2L^2_{gyy}r^4_v\delta^2_\epsilon}{\mu} \nonumber \\
 & \leq (\frac{1}{2\tau}-\frac{\mu}{4})\|v_t -v^*_t\|^2 - \frac{3}{8\tau} \|v_{t+1}-v_t\|^2 + \frac{2L^2_{gyy}r^4_v\delta^2_\epsilon}{\mu},
\end{align}
where the second inequality holds by $L_g \geq \mu$, and the last inequality is due to
$0< \tau \leq \frac{1}{6L_g}$.
It implies that
\begin{align} \label{eq:E8}
\|v_{t+1}-v^*_t\|^2 \leq ( 1 -\frac{\mu\tau}{2})\|v_t -v^*_t\|^2 - \frac{3}{4} \|v_{t+1}-v_t\|^2
+ \frac{4\tau L^2_{gyy}r^4_v\delta^2_\epsilon}{\mu}.
\end{align}

Since $v^*_t = \Big(\mathcal{S}_{[\mu,L_g]}\big[\nabla^2_{yy}g(x_t,y_t)\big]\Big)^{-1}\nabla_yf(x_t,y_t)=\arg\min_{v\in \Lambda_v}R(x_t,y_t,v)$ with $r_v=\frac{C_{fy}}{\mu}$ and \\ $v^*_{t+1} = \Big(\mathcal{S}_{[\mu,L_g]}\big[\nabla^2_{yy}g(x_{t+1},y_{t+1})\big]\Big)^{-1}\nabla_yf(x_{t+1},y_{t+1})=\arg\min_{v\in \Lambda_v}R(x_{t+1},y_{t+1},v)$, we have
\begin{align} \label{eq:E9}
 \|v^*_{t+1} -v^*_t \|^2 &= \Big\|\Big(\mathcal{S}_{[\mu,L_g]}\big[\nabla^2_{yy}g(x_{t+1},y_{t+1})\big]\Big)^{-1}\nabla_yf(x_{t+1},y_{t+1})-
 \Big(\mathcal{S}_{[\mu,L_g]}\big[\nabla^2_{yy}g(x_t,y_t)\big]\Big)^{-1}\nabla_yf(x_t,y_t)\Big\|^2 \nonumber \\
 & = \Big\| \Big(\mathcal{S}_{[\mu,L_g]}\big[\nabla^2_{yy}g(x_{t+1},y_{t+1})\big]\Big)^{-1}\nabla_yf(x_{t+1},y_{t+1})
 -\Big(\mathcal{S}_{[\mu,L_g]}\big[\nabla^2_{yy}g(x_{t+1},y_{t+1})\big]\Big)^{-1}\nabla_yf(x_t,y_t) \nonumber \\
& \quad  + \Big(\mathcal{S}_{[\mu,L_g]}\big[\nabla^2_{yy}g(x_{t+1},y_{t+1})\big]\Big)^{-1}\nabla_yf(x_t,y_t) -
 \Big(\mathcal{S}_{[\mu,L_g]}\big[\nabla^2_{yy}g(x_t,y_t)\big]\Big)^{-1}\nabla_yf(x_t,y_t)\Big\|^2 \nonumber \\
& \leq \big(\frac{4L^2_f}{\mu^2}+ \frac{4L^2_{gyy}C^2_{fx}}{\mu^4}\big)\big(\|x_{t+1}-x_t\|^2 + \|y_{t+1}-y_t\|^2\big).
\end{align}

Next, we decompose the term $\|v_{t+1} - v^*_{t+1}\|^2$ as follows:
\begin{align} \label{eq:E10}
  \|v_{t+1} - v^*_{t+1}\|^2 & = \|v_{t+1} - v^*_t + v^*_t - v^*_{t+1}\|^2   \\
  & = \|v_{t+1} - v^*_t\|^2 + 2\langle v_{t+1} - v^*_t, v^*_t - v^*_{t+1}\rangle  + \|v^*_t - v^*_{t+1}\|^2  \nonumber \\
  & \leq (1+\frac{\mu\tau}{4})\|v_{t+1} - v^*_t\|^2  + (1+\frac{4}{\mu\tau})\|v^*_t - v^*_{t+1}\|^2 \nonumber \\
  & \leq (1+\frac{\mu\tau}{4})\|v_{t+1} - v^*_t\|^2  + (1+\frac{4}{\mu\tau}) \big(\frac{4L^2_f}{\mu^2}+ \frac{4L^2_{gyy}C^2_{fx}}{\mu^4}\big)\big(\|x_{t+1}-x_t\|^2 + \|y_{t+1}-y_t\|^2\big),  \nonumber
\end{align}
where the first inequality holds by Cauchy-Schwarz inequality and Young's inequality, and the second inequality is due to the above inequality~(\ref{eq:E9}).

By combining the above inequalities \eqref{eq:E8} and \eqref{eq:E10}, we have
\begin{align}
 \|v_{t+1} - v^*_{t+1}\|^2 & \leq (1+\frac{\mu\tau}{4})( 1-\frac{\mu\tau}{2})\|v_t -v^*_t\|^2
 - (1+\frac{\mu\tau}{4})\frac{3}{4} \|v_{t+1}-v_t\|^2   \nonumber \\
 & \quad + (1+\frac{\mu\tau}{4})\frac{4\tau L^2_{gyy}r^4_v\delta^2_\epsilon}{\mu}
  + (1+\frac{4}{\mu\tau}) \big(\frac{4L^2_f}{\mu^2}+ \frac{4L^2_{gyy}C^2_{fx}}{\mu^4}\big)\big(\|x_{t+1}-x_t\|^2 + \|y_{t+1}-y_t\|^2\big).   \nonumber
\end{align}
Since $0< \tau \leq \frac{1}{6L_g}$ and $L_g \geq \mu$, we have $\tau \leq \frac{1}{6L_g} \leq \frac{1}{6\mu}$. Then we have
\begin{align}
  (1+\frac{\mu\tau}{4})(1-\frac{\mu\tau}{2})&= 1-\frac{\mu\tau}{2} +\frac{\mu\tau}{4}
  - \frac{\mu^2\tau^2}{8} \leq 1-\frac{\mu\tau}{4}, \nonumber \\
 - (1+\frac{\mu\tau}{4})\frac{3}{4} &\leq -\frac{3}{4}, \nonumber \\
  (1+\frac{\mu\tau}{4})\frac{4\tau}{\mu} & \leq (1+\frac{1}{24})\frac{4\tau}{\mu}=\frac{25\tau}{6\mu}, \nonumber \\
  1+\frac{4}{\mu\tau} & \leq \frac{5}{3}.  \nonumber
\end{align}
Thus we have
\begin{align}
    \|v_{t+1} - v^*_{t+1}\|^2 & \leq ( 1-\frac{\mu\tau}{4})\|v_t -v^*_t\|^2
 - \frac{3}{4} \|v_{t+1}-v_t\|^2   \nonumber \\
 & \quad + \frac{25\tau L^2_{gyy}r^4_v\delta^2_\epsilon}{6\mu}
  + \frac{20}{3}\big(\frac{L^2_f}{\mu^2}+ \frac{L^2_{gyy}C^2_{fx}}{\mu^4}\big)\big(\|x_{t+1}-x_t\|^2 + \|y_{t+1}-y_t\|^2\big).
\end{align}

\end{proof}

\begin{lemma}\label{lem:B2}
(Restatement of Lemma~\ref{lem:5})
Assume the sequence $\{x_t,y_t,v_t\}_{t=1}^T$ be generated from Algorithm~\ref{alg:1}, given $0<\gamma\leq \frac{1}{2L_F}$, we have
\begin{align*}
    \Phi(x_{t+1}) &\leq \Phi(x_t) - \frac{\gamma}{2}\|\mathcal{G}(x_t,w_t,\gamma)\|^2  + \frac{12\gamma}{\mu}\big(L^2_f+ r^2_vL^2_{gxy}\big)\big(g(x_t,y_t)-G(x_t)\big) \nonumber \\
    & \quad + 6\gamma C^2_{gxy}\|v_t-v^*_t\|^2 + 2\gamma L^2_{gxy}\delta^2_{\epsilon}r^4_v,
\end{align*}
where $\Phi(x) = F(x)+\phi(x)$ and $\mathcal{G}(x_t,w_t,\gamma) = \frac{1}{\gamma}(x_t-x_{t+1})$.
\end{lemma}

\begin{proof}
By the line 4 of Algorithm~\ref{alg:1}, we have
\begin{align} \label{eq:H1}
x_{t+1} = \mathbb{P}^\gamma_{\phi(\cdot)}\big(x_t, w_t\big) = \arg\min_{x\in \mathbb{R}^d}\Big\{ \langle w_t, x\rangle + \frac{1}{2\gamma}\|x-x_t\|^2 + \phi(x)\Big\}.
\end{align}

Then we define a gradient mapping $\mathcal{G}(x_t,w_t,\gamma) = \frac{1}{\gamma}(x_t-x_{t+1})$.
By the optimality condition of the subproblem~(\ref{eq:H1}), we have for any $z\in \mathbb{R}^d$
\begin{align}
 \big\langle w_t + \frac{1}{\gamma}(x_{t+1}-x_t) + \vartheta_{t+1}, z-x_{t+1}\big\rangle \geq 0,
\end{align}
where $\vartheta_{t+1}\in \partial \phi(x_{t+1})$.

Let $z=x_t$, and by the convexity of $\phi(x)$, we can obtain
\begin{align}  \label{eq:H2}
 \langle w_t, x_t - x_{t+1}\rangle & \geq \frac{1}{\gamma}\|x_{t+1}-x_t\|^2 + \langle \vartheta_{t+1}, x_{t+1}-x_t\rangle \nonumber \\
 & \geq \frac{1}{\gamma}\|x_{t+1}-x_t\|^2 + \phi(x_{t+1})-\phi(x_t).
\end{align}

According to the Lemma \ref{lem:2}, function $F(x)$ has $L_F$-Lipschitz continuous gradient.
Let $\mathcal{G}(x_t,w_t,\gamma) = \frac{1}{\gamma}(x_t-x_{t+1})$, we have
\begin{align} \label{eq:H3}
  F(x_{t+1}) & \leq F(x_t) + \langle \nabla F(x_t), x_{t+1}-x_t\rangle + \frac{L_F}{2}\|x_{t+1}-x_t\|^2 \nonumber \\
  & = F(x_t) + \langle w_t, x_{t+1}-x_t\rangle + \gamma \langle \nabla F(x_t) -w_t, \mathcal{G}(x_t,w_t,\gamma)\rangle+ \frac{\gamma ^2L_F}{2}\|\mathcal{G}(x_t,w_t,\gamma)\|^2 \nonumber \\
  & \leq F(x_t) - \gamma \|\mathcal{G}(x_t,w_t,\gamma)\|^2 - \phi(x_{t+1}) + \phi(x_t) + \gamma \langle \nabla F(x_t) - w_t, \mathcal{G}(x_t,w_t,\gamma)\rangle + \frac{\gamma ^2L_F}{2}\|\mathcal{G}(x_t,w_t,\gamma)\|^2 \nonumber \\
  & \mathop{\leq}^{(ii)} F(x_t) - \frac{\gamma}{2}\|\mathcal{G}(x_t,w_t,\gamma)\|^2 - \phi(x_{t+1}) + \phi(x_t) + \gamma\|w_t - \nabla F(x_t)\|^2,
\end{align}
where the second last inequality holds by the above inequality~(\ref{eq:H2}), and the last inequality holds by $0<\gamma\leq \frac{1}{2L_F}$ and the following inequality
\begin{align}
 \langle \nabla F(x_t) -w_t, \mathcal{G}(x_t,w_t,\gamma)\rangle &\leq \|w_t - \nabla F(x_t)\|\|\mathcal{G}(x_t,w_t,\gamma)\| \nonumber \\
 & \leq  \frac{\rho}{2}\|w_t - \nabla F(x_t)\|^2 + \frac{1}{2\rho}\|\mathcal{G}(x_t,w_t,\gamma)\|^2 \nonumber \\
 & = \|w_t - \nabla F(x_t)\|^2 + \frac{1}{4}\|\mathcal{G}(x_t,w_t,\gamma)\|^2,
\end{align}
where the above inequality holds by Young inequality with $\rho=2$.

Since $w_t = \widetilde{\nabla}f(x_t, y_t, v_t)$ in Algorithm~\ref{alg:1}, we have
\begin{align} \label{eq:H4}
\|w_t - \nabla F(x_t)\|^2 & = \|\widetilde{\nabla}f(x_t, y_t, v_t)- \nabla F(x_t)\|^2 \nonumber \\
& \leq 2\|\widehat{\nabla}f(x_t, y_t, v_t) - \nabla F(x_t)\|^2  + 2\|\widetilde{\nabla}f(x_t, y_t, v_t) - \widehat{\nabla}f(x_t, y_t, v_t) \|^2 \nonumber \\
& \leq \frac{12}{\mu}\big(L^2_f+ r^2_vL^2_{gxy}\big)\big(g(x_t,y_t)-G(x_t)\big) + 6C^2_{gxy}\|v_t-v^*_t\|^2+2L^2_{gxy}\delta^2_{\epsilon}r^4_v,
\end{align}
where the last inequality holds by the following inequalities~(\ref{eq:H5}) and~(\ref{eq:H6}).

Considering the term $\|\widehat{\nabla}f(x_t, y_t, v_t) - \nabla F(x_t)\|^2$, we have
\begin{align} \label{eq:H5}
 & \|\widehat{\nabla}f(x_t, y_t, v_t) - \nabla F(x_t)\|^2 \nonumber \\
 & = \|\nabla_xf(x_t,y_t) - \nabla^2_{xy}g(x_t,y_t)v_t - \nabla F(x_t)\|^2 \nonumber \\
 & = \|\nabla_xf(x_t,y_t) - \nabla^2_{xy}g(x_t,y_t)v_t - \nabla_xf(x_t,y^*(x_t)) + \nabla_{xy}^2g(x_t,y^*(x_t))v^*_t\|^2 \nonumber \\
 & = \|\nabla_xf(x_t,y_t) - \nabla_xf(x_t,y^*(x_t)) - \nabla^2_{xy}g(x_t,y_t)v_t + \nabla^2_{xy}g(x_t,y^*(x_t))v_t \nonumber \\
 & \quad  -\nabla^2_{xy}g(x_t,y^*(x_t))v_t+ \nabla_{xy}^2g(x_t,y^*(x_t))v^*_t\|^2 \nonumber \\
 & \leq \big(3L^2_f+ 3r^2_vL^2_{gxy}\big)\|y_t-y^*(x_t)\|^2+ 3C^2_{gxy}\|v_t-v^*_t\|^2 \nonumber \\
 & \leq \frac{6}{\mu}\big(L^2_f+ r^2_vL^2_{gxy}\big)\big(g(x_t,y_t)-G(x_t)\big) + 3C^2_{gxy}\|v_t-v^*_t\|^2,
\end{align}
where the last inequality holds by the above Lemma~\ref{lem:A2}.

Considering the term $\|\widetilde{\nabla}f(x_t, y_t, v_t) - \widehat{\nabla}f(x_t, y_t, v_t) \|^2$, we have
\begin{align} \label{eq:H6}
 & \|\widetilde{\nabla}f(x_t, y_t, v_t) - \widehat{\nabla}f(x_t, y_t, v_t) \|^2 \nonumber \\
 & = \|\nabla_xf(x_t,y_t) - \widetilde{J}(x_t, y_t, v_t, \delta_{\epsilon})-\nabla_xf(x_t,y_t) + \nabla^2_{xy}g(x_t,y_t)v_t \|^2 \nonumber \\
 & = \left\|\frac{\nabla_x g(x_t, y_t+\delta_{\epsilon} v_t) - \nabla_x g(x_t, y_t-\delta_{\epsilon} v_t)}{2\delta_{\epsilon}}- \nabla^2_{xy}g(x_t,y_t)v_t \right\|^2 \nonumber \\
 & =  \frac{1}{4\delta_{\epsilon}^2}\big\|\nabla_x g(x_t, y_t+\delta_{\epsilon} v_t) - \nabla_x g(x_t, y_t) - \delta_{\epsilon}\nabla^2_{xy}g(x_t,y_t)v_t \nonumber \\
 & \qquad + \nabla_x g(x_t, y_t)- \nabla_x g(x_t, y_t-\delta_{\epsilon} v_t)- \delta_{\epsilon}\nabla^2_{xy}g(x_t,y_t)v_t \big\|^2 \nonumber \\
 & \leq \frac{1}{2\delta_{\epsilon}^2}\big\|\nabla_x g(x_t, y_t+\delta_{\epsilon} v_t) - \nabla_x g(x_t, y_t) - \delta_{\epsilon}\nabla^2_{xy}g(x_t,y_t)v_t \big\|^2 \nonumber \\
 & \quad + \frac{1}{2\delta_{\epsilon}^2}\big\|\nabla_x g(x_t, y_t)- \nabla_x g(x_t, y_t-\delta_{\epsilon} v_t)- \delta_{\epsilon}\nabla^2_{xy}g(x_t,y_t)v_t \big\|^2 \nonumber \\
 & = \frac{1}{2\delta_{\epsilon}^2}\Big\|\int_{k=0}^1\big(\nabla^2_{xy} g(x_t, y_t+k\delta_{\epsilon} v_t) -\nabla^2_{xy}g(x_t,y_t)\big)\delta_{\epsilon}v_t dk \Big\|^2 \nonumber \\
 & \quad + \frac{1}{2\delta_{\epsilon}^2}\Big\|\int_{k=0}^1\big(\nabla^2_{xy} g(x_t, y_t-k\delta_{\epsilon} v_t) -\nabla^2_{xy}g(x_t,y_t)\big)\delta_{\epsilon}v_t dk  \Big\|^2 \nonumber \\
 & \leq \frac{1}{2\delta_{\epsilon}^2}\Big(\int_{k=0}^1\big\|\nabla^2_{xy} g(x_t, y_t+k\delta_{\epsilon} v_t) -\nabla^2_{xy}g(x_t,y_t)\big\|\|v_t\| \delta_{\epsilon}dk\Big)^2 \nonumber \\
 & \quad + \frac{1}{2\delta_{\epsilon}^2}\Big(\int_{k=0}^1\big\|\nabla^2_{xy} g(x_t, y_t-k\delta_{\epsilon} v_t) -\nabla^2_{xy}g(x_t,y_t)\big\|\|v_t\|\delta_{\epsilon}dk \Big)^2 \nonumber \\
 & \leq L^2_{gxy}\delta^2_{\epsilon}\|v_t\|^4\leq L^2_{gxy}\delta^2_{\epsilon}r^4_v,
\end{align}
where the last inequality holds by $\|v_t\| \leq r_v$.

By combining the above inequalities~(\ref{eq:H3}) with~(\ref{eq:H4}), we can obtain
\begin{align*}
    \Phi(x_{t+1}) &\leq \Phi(x_t) - \frac{\gamma}{2}\|\mathcal{G}(x_t,w_t,\gamma)\|^2  + \frac{12\gamma}{\mu}\big(L^2_f+ r^2_vL^2_{gxy}\big)\big(g(x_t,y_t)-G(x_t)\big) \nonumber \\
    & \quad + 6\gamma C^2_{gxy}\|v_t-v^*_t\|^2 + 2\gamma L^2_{gxy}\delta^2_{\epsilon}r^4_v.
\end{align*}

\end{proof}

\begin{lemma} \label{lem:B3}
(Restatement of Lemma~\ref{lem:6})
Suppose the sequence $\{x_t,y_t,v_t\}_{t=1}^T$ be generated from Algorithm~\ref{alg:1}.
Under the above Assumptions~\ref{ass:1}-\ref{ass:3}, given $\gamma\leq \min\Big\{\frac{\lambda\mu }{16L_G }, \frac{\mu }{16L^2_g } \Big\}$ and $0<\lambda \leq \frac{1}{2L_g }$, we have
\begin{align}
g(x_{t+1},y_{t+1}) - G(x_{t+1})
& \leq (1-\frac{\lambda\mu}{2 }) \big(g(x_t,y_t) -G(x_t)\big) + \frac{1}{8\gamma}\|x_{t+1}-x_t\|^2  -\frac{1}{4\lambda }\|y_{t+1}-y_t\|^2 \nonumber \\
& \qquad + \lambda\|\nabla_y g(x_t,y_t)-u_t\|^2,
\end{align}
where $G(x_t)=g(x_t,y^*(x_t))$ with $y^*(x_t) \in \arg\min_{y}g(x_t,y)$ for all $t\geq 1$.
\end{lemma}

\begin{proof}
Using the Assumption~\ref{ass:3}, i.e., $L_g$-smoothness of $g(x,\cdot)$, such that
\begin{align} \label{eq:EE1}
    g(x_{t+1},y_{t+1}) \leq g(x_{t+1},y_t) + \langle \nabla_y g(x_{t+1},y_t), y_{t+1}-y_t \rangle + \frac{L_g}{2}\|y_{t+1}-y_t\|^2.
\end{align}

Since $y_{t+1} = \arg\min_{y\in \mathbb{R}^p}\big\{ \langle u_t, y\rangle
+ \frac{1}{2\lambda}(y-y_t)^T(y-y_t) \big\} = y_t - \lambda u_t$, we can obtain
\begin{align} \label{eq:EE2}
    & \langle \nabla_y g(x_{t+1},y_t), y_{t+1}-y_t \rangle \nonumber \\
    & = -\lambda\langle \nabla_y g(x_{t+1},y_t), u_t \rangle \nonumber \\
    & = -\frac{\lambda}{2}\Big( \|\nabla_y g(x_{t+1},y_t)\|^2 + \|u_t\|^2 - \|\nabla_y g(x_{t+1},y_t)-\nabla_y g(x_t,y_t) + \nabla_y g(x_t,y_t)-u_t\|^2 \Big) \nonumber \\
    & \leq -\frac{\lambda}{2 } \|\nabla_y g(x_{t+1},y_t)\|^2 -\frac{1}{2\lambda } \|y_{t+1}-y_t\|^2 + \lambda L^2_g\|x_{t+1}-x_t\|^2 + \lambda\|\nabla_y g(x_t,y_t)-u_t\|^2 \nonumber \\
    & \leq -\lambda\mu\big(g(x_{t+1},y_t)-G(x_{t+1})\big)-\frac{1}{2\lambda} \|y_{t+1}-y_t\|^2 + \lambda L^2_g \|x_{t+1}-x_t\|^2 + \lambda\|\nabla_y g(x_t,y_t)-u_t\|^2,
\end{align}
where the last inequality is due to the quadratic growth condition of $\mu$-PL functions, i.e.,
\begin{align}
    \|\nabla_y g(x_{t+1},y_t)\|^2 \geq 2\mu\big( g(x_{t+1},y_t)-\min_{y'}g(x_{t+1},y')\big) = 2\mu\big( g(x_{t+1},y_t) - G(x_{t+1})\big).
\end{align}
Substituting \eqref{eq:EE2} into \eqref{eq:EE1}, we have
\begin{align} \label{eq:EE3}
    g(x_{t+1},y_{t+1})
    & \leq g(x_{t+1},y_t)-\lambda\mu\big(g(x_{t+1},y_t)- G(x_{t+1})\big)-\frac{1}{2\lambda} \|y_{t+1}-y_t\|^2 + \lambda L^2_g\|x_{t+1}-x_t\|^2 \nonumber \\
    & \quad + \lambda\|\nabla_y g(x_t,y_t)-u_t\|^2 + \frac{L_g}{2}\|y_{t+1}-y_t\|^2,
\end{align}
then rearranging the terms, we can obtain
\begin{align} \label{eq:EE4}
    g(x_{t+1},y_{t+1})- G(x_{t+1})
    & \leq (1-\lambda\mu)\big(g(x_{t+1},y_t)- G(x_{t+1})\big)-\frac{1}{2\lambda } \|y_{t+1}-y_t\|^2 + \lambda L^2_g \|x_{t+1}-x_t\|^2 \nonumber \\
    & \quad + \lambda\|\nabla_y g(x_t,y_t)-u_t\|^2 + \frac{L_g}{2}\|y_{t+1}-y_t\|^2.
\end{align}

Next, using $L_g$-smoothness of function $f(\cdot,y_t)$, such that
\begin{align}
     g(x_{t+1},y_t) \leq g(x_t,y_t) + \langle \nabla_x g(x_t,y_t), x_{t+1}-x_t \rangle +  \frac{L_g}{2}\|x_{t+1}-x_t\|^2 ,
\end{align}
then we have
\begin{align}
    &g(x_{t+1},y_t) - g(x_t,y_t) \nonumber \\
    & \leq \langle \nabla_x g(x_t,y_t), x_{t+1}-x_t \rangle + \frac{L_g}{2}\|x_{t+1}-x_t\|^2 \nonumber \\
    & = \langle \nabla_x g(x_t,y_t) - \nabla G(x_t), x_{t+1}-x_t \rangle + \langle \nabla G(x_t), x_{t+1}-x_t \rangle + \frac{L_g}{2}\|x_{t+1}-x_t\|^2 \nonumber \\
    & \leq \frac{1}{8\gamma}\|x_{t+1}-x_t\|^2 + 2\gamma\|\nabla_x g(x_t,y_t) - \nabla G(x_t)\|^2  + \langle \nabla G(x_t), x_{t+1}-x_t \rangle + \frac{L_g}{2}\|x_{t+1}-x_t\|^2 \nonumber \\
    & \leq \frac{1}{8\gamma}\|x_{t+1}-x_t\|^2 + 2L^2_g\gamma \|y_t - y^*(x_t)\|^2 + G(x_{t+1}) - G(x_t) + \frac{L_G}{2}\|x_{t+1}-x_t\|^2 + \frac{L_g}{2}\|x_{t+1}-x_t\|^2 \nonumber \\
    & \leq \frac{4L^2_g\gamma}{\mu } \big(g(x_t,y_t) - G(x_t)\big) + G(x_{t+1})- G(x_t)+ (\frac{1 }{8\gamma}+L_G)\|x_{t+1}-x_t\|^2,
\end{align}
where the second last inequality is due to
$L_G$-smoothness of function $G(x)$, and the last inequality holds by Lemma~\ref{lem:A2} and $L_g\leq L_G$.
Then we have
\begin{align} \label{eq:EE5}
    g(x_{t+1},y_t) -G(x_{t+1}) & = g(x_{t+1},y_t)- g(x_t,y_t) + g(x_t,y_t)- G(x_t) + G(x_t) -G(x_{t+1})
    \nonumber \\
    & \leq (1+\frac{4L^2_g\gamma}{\mu }) \big(g(x_t,y_t) -G(x_t)\big) + (\frac{1}{8\gamma}+ L_G)\|x_{t+1}-x_t\|^2.
\end{align}

Substituting \eqref{eq:EE5} in \eqref{eq:EE4}, we get
\begin{align}
    & g(x_{t+1},y_{t+1})- G(x_{t+1})\nonumber \\
    & \leq (1-\lambda\mu)(1+\frac{4L^2_g\gamma}{\mu}) \big(g(x_t,y_t) -G(x_t)\big) + (1-\lambda\mu)(\frac{ 1}{8\gamma}+L_G)\|x_{t+1}-x_t\|^2 \nonumber \\
    & \quad -\frac{1}{2\lambda} \|y_{t+1}-y_t\|^2 + \lambda L^2_g\|x_{t+1}-x_t\|^2 + \lambda\|\nabla_y g(x_t,y_t)-u_t\|^2 + \frac{L_g}{2}\|y_{t+1}-y_t\|^2 \nonumber \\
    & = (1-\lambda\mu)(1+\frac{4L^2_g\gamma}{\mu}) \big(g(x_t,y_t) -G(x_t)\big) + \big(\frac{ 1}{8\gamma}+L_G-\frac{\lambda\mu}{8\gamma }-L_G\lambda\mu +L^2_g\lambda\big)\|x_{t+1}-x_t\|^2 \nonumber \\
    & \quad -\frac{1}{2}\big(\frac{1}{\lambda }-L_g\big) \|y_{t+1}-y_t\|^2 + \lambda\|\nabla_y g(x_t,y_t)-u_t\|^2 \nonumber \\
    & \leq (1-\frac{\lambda\mu}{2}) \big(g(x_t,y_t) -G(x_t)\big) + \frac{1}{8\gamma}\|x_{t+1}-x_t\|^2  -\frac{1}{4\lambda}\|y_{t+1}-y_t\|^2 + \lambda\|\nabla_y g(x_t,y_t)-u_t\|^2,
\end{align}
where the last inequality holds by $\gamma\leq \min\Big\{\frac{\lambda\mu }{16L_G}, \frac{\mu }{16L^2_g}\Big\}$, $L_G\geq L_g(1+\kappa)^2$ and $\lambda\leq \frac{1}{2L_g}$ for all $t\geq 1$, i.e.,
\begin{align}
   & \gamma\leq \frac{\lambda\mu }{16L_G } \Rightarrow \lambda \geq \frac{16L_G\gamma }{\mu } \geq   \frac{16L_g}{\mu }(1+\kappa)^2\gamma  \geq 8\kappa^2\gamma \Rightarrow  \frac{\lambda\mu}{2 } \geq \frac{4L^2_g\gamma}{\mu } \nonumber \\
   & \gamma \leq \min\Big\{\frac{\lambda\mu }{16L_G }, \frac{\mu }{16L^2_g } \Big\}
   \Rightarrow \frac{\lambda\mu}{8\gamma} \geq
    L_G+L^2_g\lambda \nonumber \\
   &\lambda \leq \frac{1}{2L_g }  \Rightarrow \frac{1}{2\lambda } \geq  L_g, \ \forall t\geq 1.
\end{align}

\end{proof}

\begin{theorem}  \label{th:A1}
(Restatement of Theorem~\ref{th:1})
 Assume the sequence $\{x_t,y_t,v_t\}_{t=1}^T$ be generated from our Algorithm \ref{alg:1}. Under the above Assumptions~\ref{ass:1}-\ref{ass:5}, let $0<\gamma\leq \min\Big(\frac{1}{2L_F},\frac{\lambda\mu }{16L_G}, \frac{\mu }{16L^2_g},\frac{3}{160\breve{L}^2}, \frac{\mu\tau}{30C^2_{gxy}},\frac{\mu^2\lambda}{30(L^2_f+ r^2_vL^2_{gxy})}\Big)$, $0< \lambda\leq \min\big(\frac{1}{2L_g},\frac{3}{80\breve{L}^2}\big)$ and $0<\tau \leq \frac{1}{6L_g}$,
 we have
 \begin{align}
 & \frac{1}{T}\sum_{t=1}^T\|\mathcal{G}(x_t,\nabla F(x_t),\gamma)\|^2 \nonumber \\
 & \leq \frac{8(\Phi(x_1) +  g(x_1,y_1)-G(x_1)+ \|v_1-v_1^*\|^2 - \Phi^*)}{T\gamma} + 20 L^2_{gxy}\delta^2_{\epsilon}r^4_v + \frac{100\tau L^2_{gyy}r^4_v\delta^2_\epsilon}{3\gamma\mu},
\end{align}
where $\Phi(x)=F(x)+\phi(x)$ and $\breve{L}^2=\frac{L^2_f}{\mu^2}+ \frac{L^2_{gyy}C^2_{fx}}{\mu^4}$.
\end{theorem}

\begin{proof}
According to Lemma \ref{lem:B2}, we have
\begin{align} \label{eq:D1}
    \Phi(x_{t+1}) - \Phi(x_t) & \leq - \frac{\gamma}{2}\|\mathcal{G}(x_t,w_t,\gamma)\|^2  + \frac{12\gamma}{\mu}\big(L^2_f+ r^2_vL^2_{gxy}\big)\big(g(x_t,y_t)-G(x_t)\big)  \nonumber \\
    & \quad + 6\gamma C^2_{gxy}\|v_t-v^*_t\|^2+ 2\gamma L^2_{gxy}\delta^2_{\epsilon}r^4_v,
\end{align}

According to the line 3 of Algorithm~\ref{alg:1}, we have $u_t=\nabla_y g(x_t,y_t)$.
Then by using Lemma \ref{lem:B3}, we have
\begin{align} \label{eq:D2}
& g(x_{t+1},y_{t+1}) - G(x_{t+1}) - \big(g(x_t,y_t) -G(x_t)\big) \nonumber \\
& \leq  -\frac{\lambda\mu}{2} \big(g(x_t,y_t) -G(x_t)\big) + \frac{1}{8\gamma}\|x_{t+1}-x_t\|^2  -\frac{1}{4\lambda }\|y_{t+1}-y_t\|^2 + \lambda\|\nabla_y g(x_t,y_t)-u_t\|^2 \nonumber \\
& = -\frac{\lambda\mu}{2} \big(g(x_t,y_t) -G(x_t)\big) + \frac{\gamma}{8}\|\mathcal{G}(x_t,w_t,\gamma)\|^2
-\frac{1}{4\lambda }\|y_{t+1}-y_t\|^2,
\end{align}
where the above equality holds by $u_t=\nabla_y g(x_t,y_t)$ and $\mathcal{G}(x_t,w_t,\gamma) = \frac{1}{\gamma}(x_t-x_{t+1})$.

By using Lemma \ref{lem:B1}, we have
\begin{align} \label{eq:D3}
& \|v_{t+1} - v^*_{t+1}\|^2 -\|v_t -v^*_t\|^2 \nonumber \\
& \leq  -\frac{\mu\tau}{4}\|v_t -v^*_t\|^2
 - \frac{3}{4} \|v_{t+1}-v_t\|^2
 + \frac{25\tau L^2_{gyy}r^4_v\delta^2_\epsilon}{6\mu}\nonumber \\
 & \quad + \frac{20}{3}\big(\frac{L^2_f}{\mu^2}+ \frac{L^2_{gyy}C^2_{fx}}{\mu^4}\big)\big(\|x_{t+1}-x_t\|^2 + \|y_{t+1}-y_t\|^2\big) \nonumber \\
 &=  -\frac{\mu\tau}{4}\|v_t -v^*_t\|^2
 - \frac{3}{4} \|v_{t+1}-v_t\|^2
 + \frac{25\tau L^2_{gyy}r^4_v\delta^2_\epsilon}{6\mu} + \frac{20}{3}\breve{L}^2\gamma^2\|\mathcal{G}(x_t,w_t,\gamma)\|^2 + \frac{20}{3}\breve{L}^2\|y_{t+1}-y_t\|^2,
\end{align}
where the above equality holds by $\breve{L}^2=\frac{L^2_f}{\mu^2}+ \frac{L^2_{gyy}C^2_{fx}}{\mu^4}$ and $\mathcal{G}(x_t,w_t,\gamma) = \frac{1}{\gamma}(x_t-x_{t+1})$.

Next we define a useful Lyapunov function (i.e. potential function) for any $t\geq 1$
\begin{align}
 \Psi_t = \Phi(x_t) + g(x_t,y_t)-G(x_t)+ \|v_t-v_t^*\|^2.
\end{align}

By using the above inequalities~(\ref{eq:D1}),~(\ref{eq:D2}) and~(\ref{eq:D3}), we have
\begin{align} \label{eq:D4}
 \Psi_{t+1} - \Psi_t & = \Phi(x_{t+1}) - \Phi(x_t) + g(x_{t+1},y_{t+1})-G(x_{t+1})-\big(g(x_t,y_t)-G(x_t)\big) + \|v_{t+1}-v_{t+1}^*\|^2 - \|v_t-v_t^*\|^2
 \nonumber \\
 & \leq- \frac{\gamma}{2}\|\mathcal{G}(x_t,w_t,\gamma)\|^2  + \frac{12\gamma}{\mu}\big(L^2_f+ r^2_vL^2_{gxy}\big)\big(g(x_t,y_t)-G(x_t)\big) + 6\gamma C^2_{gxy}\|v_t-v^*_t\|^2 + 2\gamma L^2_{gxy}\delta^2_{\epsilon}r^4_v   \nonumber \\
 & \quad -\frac{\lambda\mu}{2} \big(g(x_t,y_t) -G(x_t)\big) + \frac{\gamma}{8}\|\mathcal{G}(x_t,w_t,\gamma)\|^2
-\frac{1}{4\lambda }\|y_{t+1}-y_t\|^2 \nonumber \\
 & \quad -\frac{\mu\tau}{4}\|v_t -v^*_t\|^2
 - \frac{3}{4} \|v_{t+1}-v_t\|^2 + \frac{25\tau L^2_{gyy}r^4_v\delta^2_\epsilon}{6\mu}
 + \frac{20}{3}\breve{L}^2\gamma^2\|\mathcal{G}(x_t,w_t,\gamma)\|^2 + \frac{20}{3}\breve{L}^2\|y_{t+1}-y_t\|^2 \nonumber \\
& \leq -\big(\frac{\lambda\mu}{2}- \frac{12\gamma}{\mu}\big(L^2_f+ r^2_vL^2_{gxy}\big)\big)\big(g(x_t,y_t) -G(x_t)\big) - \frac{\gamma}{4}\|\mathcal{G}(x_t,w_t,\gamma)\|^2 - \big(\frac{\mu\tau}{4}-6\gamma C^2_{gxy}\big)\|v_t -v^*_t\|^2 \nonumber \\
& \quad + 2\gamma L^2_{gxy}\delta^2_{\epsilon}r^4_v + \frac{25\tau L^2_{gyy}r^4_v\delta^2_\epsilon}{6\mu} ,
\end{align}
where the last inequality is due to $0<\gamma\leq \frac{3}{160\breve{L}^2}$ and $0<\lambda\leq \frac{3}{80\breve{L}^2}$.

Let
\begin{align}
& x^+_{t+1} = \mathbb{P}_{\phi(\cdot)}^{\gamma}(x_t, \nabla F(x_t))= \arg\min_{x\in \mathbb{R}^d}\Big\{ \langle \nabla F(x_t), x\rangle
+ \frac{1}{2\gamma}\|x-x_t\|^2 + \phi(x)\Big\}, \label{eq:D5} \\
& x_{t+1} = \mathbb{P}_{\phi(\cdot)}^{\gamma}(x_t, w_t)= \arg\min_{x\in \mathbb{R}^d}\Big\{ \langle w_t, x\rangle
+ \frac{1}{2\gamma}\|x-x_t\|^2 + \phi(x)\Big\}. \label{eq:D6}
\end{align}
By the optimality conditions of (\ref{eq:D5}) and (\ref{eq:D6}), for any $z\in \mathbb{R}^d$, there exist $\vartheta_1\in \partial\phi(x^+_{t+1})$
and $\vartheta_2\in \partial\phi(x_{t+1})$ such that
\begin{align}
& \langle \nabla F(x_t) + \frac{1}{\gamma}(x^+_{t+1}-x_t) + \vartheta_1, z- x^+_{t+1}\rangle \geq 0, \label{eq:D7} \\
& \langle w_t + \frac{1}{\gamma}(x_{t+1}-x_t) + \vartheta_2, z- x_{t+1}\rangle \geq 0. \label{eq:D8}
\end{align}
Putting $z=x_{t+1}$ into (\ref{eq:D7}), by the convexity of $\phi(x)$, we have
\begin{align} \label{eq:D9}
 \langle \nabla F(x_t), x_{t+1}- x^+_{t+1}\rangle & \geq  \frac{1}{\gamma}\langle x^+_{t+1}-x_t, x^+_{t+1}- x_{t+1}\rangle + \langle \vartheta_1, x^+_{t+1}- x_{t+1}\rangle\nonumber \\
 & \geq \frac{1}{\gamma}\langle x^+_{t+1}-x_t, x^+_{t+1}- x_{t+1}\rangle + \phi(x^+_{t+1}) - \phi(x_{t+1}).
\end{align}
Similarly, putting $z=x^+_{t+1}$ into (\ref{eq:D8}), by the convexity of $\phi(x)$, we have
\begin{align} \label{eq:D10}
 \langle w_t, x^+_{t+1}-x_{t+1} \rangle & \geq  \frac{1}{\gamma}\langle x_{t+1}-x_t, x_{t+1}- x^+_{t+1}\rangle + \langle \vartheta_2, x_{t+1}- x^+_{t+1}\rangle\nonumber \\
 & \geq \frac{1}{\gamma}\langle x_{t+1}-x_t, x_{t+1}- x^+_{t+1}\rangle +  \phi(x_{t+1})-\phi(x^+_{t+1}).
\end{align}
Summing up (\ref{eq:D9}) and (\ref{eq:D10}),  we can obtain
\begin{align}
   \|\nabla F(x_t) -w_t\|\|x_{t+1}- x^+_{t+1}\| \geq \langle \nabla F(x_t) -w_t, x_{t+1}- x^+_{t+1}\rangle  \geq  \frac{1}{\gamma}\|x^+_{t+1}- x_{t+1}\|^2.
\end{align}
Then we have
\begin{align} \label{eq:D11}
   \|\nabla F(x_t) -w_t\|  \geq  \frac{1}{\gamma}\|x^+_{t+1}- x_{t+1}\|=\frac{1}{\gamma}\left\|\mathbb{P}_{\phi(\cdot)}^{\gamma}(x_t,\nabla F(x_t))- \mathbb{P}_{\phi(\cdot)}^{\gamma}(x_t, w_t)\right\|.
\end{align}

Since $\mathcal{G}(x_t,w_t,\gamma) = \frac{1}{\gamma}\big(x_t-\mathbb{P}_{\phi(\cdot)}^{\gamma}(x_t, w_t)\big)$ and $\mathcal{G}(x_t,\nabla F(x_t),\gamma) = \frac{1}{\gamma}\big(x_t-\mathbb{P}_{\phi(\cdot)}^{\gamma}(x_t,\nabla F(x_t))\big)$,
we have
\begin{align}
 \|\mathcal{G}(x_t,\nabla F(x_t),\gamma)\|^2 & \leq 2\|\mathcal{G}(x_t,w_t,\gamma)\|^2 + 2\|\mathcal{G}(x_t,w_t,\gamma)-\mathcal{G}(x_t,\nabla F(x_t),\gamma)\|^2 \nonumber \\
 & = 2\|\mathcal{G}(x_t,w_t,\gamma)\|^2 + \frac{2}{\gamma^2}\|\mathbb{P}_{\phi(\cdot)}^{\gamma}(x_t,\nabla F(x_t))-\mathbb{P}_{\phi(\cdot)}^{\gamma}(x_t, w_t)\|^2 \nonumber \\
 & \mathop{\leq}^{(i)} 2\|\mathcal{G}(x_t,w_t,\gamma)\|^2 + 2\|w_t - \nabla F(x_t))\|^2 \nonumber \\
 & \leq 2\|\mathcal{G}(x_t,w_t,\gamma)\|^2 + \frac{24}{\mu}\big(L^2_f+ r^2_vL^2_{gxy}\big)\big(g(x_t,y_t)-G(x_t)\big) + 12C^2_{gxy}\|v_t-v^*_t\|^2 \nonumber \\
 & \quad +4L^2_{gxy}\delta^2_{\epsilon}r^4_v,
\end{align}
where the inequality (i) holds by the above inequality~(\ref{eq:D11}).
Then we can obtain
\begin{align} \label{eq:D12}
-\|\mathcal{G}(x_t,w_t,\gamma)\|^2 & \leq
-\frac{1}{2}\|\mathcal{G}(x_t,\nabla F(x_t),\gamma)\|^2 + \frac{12}{\mu}\big(L^2_f+ r^2_vL^2_{gxy}\big)\big(g(x_t,y_t)-G(x_t)\big) + 6C^2_{gxy}\|v_t-v^*_t\|^2 \nonumber \\
 & \quad +2L^2_{gxy}\delta^2_{\epsilon}r^4_v.
\end{align}

Plugging the above inequalities~(\ref{eq:D12}) into (\ref{eq:D4}),
we can further get
\begin{align} \label{eq:D13}
& \Psi_{t+1} - \Psi_t \nonumber \\
& \leq -\big(\frac{\lambda\mu}{2}- \frac{12\gamma}{\mu}\big(L^2_f+ r^2_vL^2_{gxy}\big)\big)\big(g(x_t,y_t) -G(x_t)\big) - \frac{\gamma}{4}\|\mathcal{G}(x_t,w_t,\gamma)\|^2 - \big(\frac{\mu\tau}{4}-6\gamma C^2_{gxy}\big)\|v_t -v^*_t\|^2 \nonumber \\
& \quad + 2\gamma L^2_{gxy}\delta^2_{\epsilon}r^4_v + \frac{25\tau L^2_{gyy}r^4_v\delta^2_\epsilon}{6\mu} \nonumber \\
& \leq -\big(\frac{\lambda\mu}{2}- \frac{12\gamma}{\mu}\big(L^2_f+ r^2_vL^2_{gxy}\big)\big)\big(g(x_t,y_t) -G(x_t)\big) -\frac{\gamma}{8}\|\mathcal{G}(x_t,\nabla F(x_t),\gamma)\|^2  \nonumber \\
 & \quad + \frac{3\gamma}{\mu}\big(L^2_f+ r^2_vL^2_{gxy}\big)\big(g(x_t,y_t)-G(x_t)\big) + \frac{3\gamma C^2_{gxy}}{2}\|v_t-v^*_t\|^2 + \frac{\gamma L^2_{gxy}\delta^2_{\epsilon}r^4_v}{2}  \nonumber \\
& \quad - \big(\frac{\mu\tau}{4}-6\gamma C^2_{gxy}\big)\|v_t -v^*_t\|^2 + 2\gamma L^2_{gxy}\delta^2_{\epsilon}r^4_v + \frac{25\tau L^2_{gyy}r^4_v\delta^2_\epsilon}{6\mu} \nonumber \\
& = -\big(\frac{\lambda\mu}{2}- \frac{15\gamma}{\mu}\big(L^2_f+ r^2_vL^2_{gxy}\big)\big)\big(g(x_t,y_t) -G(x_t)\big) -\frac{\gamma}{8}\|\mathcal{G}(x_t,\nabla F(x_t),\gamma)\|^2  \nonumber \\
& \quad - \big(\frac{\mu\tau}{4}-\frac{15\gamma C^2_{gxy}}{2}\big)\|v_t -v^*_t\|^2 + \frac{5\gamma L^2_{gxy}\delta^2_{\epsilon}r^4_v}{2} + \frac{25\tau L^2_{gyy}r^4_v\delta^2_\epsilon}{6\mu} \nonumber \\
& \leq -\frac{\gamma}{8}\|\mathcal{G}(x_t,\nabla F(x_t),\gamma)\|^2 + \frac{5\gamma L^2_{gxy}\delta^2_{\epsilon}r^4_v}{2} + \frac{25\tau L^2_{gyy}r^4_v\delta^2_\epsilon}{6\mu},
\end{align}
where the last inequality holds by $0<\gamma\leq \big(\frac{\mu\tau}{30C^2_{gxy}},\frac{\mu^2\lambda}{30(L^2_f+ r^2_vL^2_{gxy})}\big)$.

Based on the inequality~(\ref{eq:D13}), we have
\begin{align}
\frac{1}{T}\sum_{t=1}^T\|\mathcal{G}(x_t,\nabla F(x_t),\gamma)\|^2 & \leq \frac{1}{T}\sum_{t=1}^T\frac{8(\Psi_t  -\Psi_{t+1})}{\gamma} + 20 L^2_{gxy}\delta^2_{\epsilon}r^4_v + \frac{100\tau L^2_{gyy}r^4_v\delta^2_\epsilon}{3\gamma\mu} \nonumber \\
& \mathop{\leq}^{(i)} \frac{8(\Psi_1  - \Phi^*)}{T\gamma} + 20L^2_{gxy}\delta^2_{\epsilon}r^4_v + \frac{100\tau L^2_{gyy}r^4_v\delta^2_\epsilon}{3\gamma\mu} \nonumber \\
& = \frac{8(\Phi(x_1) +  g(x_1,y_1)-G(x_1)+ \|v_1-v_1^*\|^2 - \Phi^*)}{T\gamma} + 20L^2_{gxy}\delta^2_{\epsilon}r^4_v \nonumber \\
& \quad + \frac{100\tau L^2_{gyy}r^4_v\delta^2_\epsilon}{3\gamma\mu},
\end{align}
where the above inequality (i) holds by Assumption~\ref{ass:5}.

Set $\delta_{\epsilon}=O\big(\frac{1}{\sqrt{T}\max(L^2_{gxy},L^2_{gyy}/\mu)r^2_v}\big)$, we can obtain
\begin{align}
\min_{1\leq t\leq T}\|\mathcal{G}(x_t,\nabla F(x_t),\gamma)\|^2 \leq \frac{1}{T}\sum_{t=1}^T\|\mathcal{G}(x_t,\nabla F(x_t),\gamma)\|^2 & \leq O(\frac{1}{T}).
\end{align}

\end{proof}

\subsection{ Convergence Analysis of of HJFBiO Algorithm for Bilevel Optimization without Regularization}
\label{appendix:C2}
In this subsection, we provide the convergence analysis of our HJFBiO algorithm for bilevel optimization without Regularization.

\begin{lemma}\label{lem:B4}
Assume the sequence $\{x_t,y_t,v_t\}_{t=1}^T$ be generated from Algorithm~\ref{alg:1}, given $0<\gamma\leq \frac{1}{2L_F}$, we have
\begin{align*}
    F(x_{t+1}) &\leq F(x_{t}) - \frac{\gamma}{2}\|\nabla F(x_t)\|^2 - \frac{\gamma}{4}\|\widetilde{\nabla}f(x_t, y_t, v_t)\|^2 + \gamma L^2_{gxy}r^4_v\delta^2_{\epsilon} \nonumber \\
   & \quad + \frac{6\gamma}{\mu}\big(L^2_f+ r^2_vL^2_{gxy}\big)\big(g(x_t,y_t)-G(x_t)\big)+ 3\gamma C^2_{gxy}\|v_t-v^*_t\|^2,
\end{align*}

\end{lemma}

\begin{proof}
When $\phi(x)\equiv 0$, at the line 4 of Algorithm~\ref{alg:1}, we have $x_{t+1}=x_t-\gamma w_t$.
By using the Lipschitz smoothness of the objective function $F(x)$, we have
\begin{align}\label{eq:BB1}
   F(x_{t+1}) &\leq F(x_{t}) + \langle \nabla F(x_t), x_{t+1} - x_t\rangle + \frac{L_F}{2}\|x_{t+1} - x_t\|^2 \nonumber \\
   & = F(x_t) -\gamma \langle \nabla F(x_t), w_t\rangle + \frac{\gamma^2L_F}{2}\|w_t\|^2 \nonumber \\
   & = F(x_{t}) - \frac{\gamma}{2}\|\nabla F(x_t)\|^2 - \frac{\gamma}{2}(1-\gamma L_F)\|w_t\|^2 + \frac{\gamma}{2}\|w_t - \nabla F(x_t)\|^2 \nonumber \\
   & = F(x_{t}) - \frac{\gamma}{2}\|\nabla F(x_t)\|^2 - \frac{\gamma}{2}(1-\gamma L_F)\|\widetilde{\nabla}f(x_t, y_t, v_t)\|^2 + \frac{\gamma}{2}\|\widetilde{\nabla}f(x_t, y_t, v_t) - \nabla F(x_t)\|^2 \nonumber \\
   & \mathop{\leq}^{(i)} F(x_{t}) - \frac{\gamma}{2}\|\nabla F(x_t)\|^2 - \frac{\gamma}{4}\|\widetilde{\nabla}f(x_t, y_t, v_t)\|^2 + \frac{\gamma}{2}\|\widetilde{\nabla}f(x_t, y_t, v_t) - \nabla F(x_t)\|^2 \nonumber \\
   & \mathop{\leq}^{(ii)} F(x_{t}) - \frac{\gamma}{2}\|\nabla F(x_t)\|^2 - \frac{\gamma}{4}\|\widetilde{\nabla}f(x_t, y_t, v_t)\|^2 + \gamma L^2_{gxy}r^4_v\delta^2_{\epsilon} \nonumber \\
   & \quad + \frac{6\gamma}{\mu}\big(L^2_f+ r^2_vL^2_{gxy}\big)\big(g(x_t,y_t)-G(x_t)\big)+ 3\gamma C^2_{gxy}\|v_t-v^*_t\|^2,
\end{align}
where the above inequality (i) holds by $0<\gamma\leq \frac{1}{2L_F}$, and the above inequality (ii) holds by the above inequality~(\ref{eq:H4}).

\end{proof}

\begin{theorem}  \label{th:A2}
 Assume the sequence $\{x_t,y_t,v_t\}_{t=1}^T$ be generated from our Algorithm \ref{alg:1}. Under the above Assumptions~\ref{ass:1}-\ref{ass:5}, let $0<\gamma\leq \min\Big(\frac{1}{2L_F},\frac{\lambda\mu }{16L_G}, \frac{\mu }{16L^2_g},\frac{3}{160\breve{L}^2}, \frac{\mu\tau}{12C^2_{gxy}},\frac{\lambda\mu^2}{12\big(L^2_f+ r^2_vL^2_{gxy}\big)}\Big)$, $0< \lambda\leq \min\big(\frac{1}{2L_g},\frac{3}{80\breve{L}^2}\big)$ and $0<\tau \leq \frac{1}{6L_g}$,
 we have
 \begin{align}
& \frac{1}{T}\sum_{t=1}^T\|\nabla F(x_t)\|^2  \nonumber \\
& \leq \frac{2\big(F(x_1) + g(x_1,y_1)-G(x_1) + \|v_1-v_1^*\|^2  - F^*\big)}{T\gamma} + 2L^2_{gxy}r^4_v\delta^2_{\epsilon}+ \frac{25\tau L^2_{gyy}r^4_v\delta^2_\epsilon}{2\gamma\mu},
\end{align}
where $\breve{L}^2=\frac{L^2_f}{\mu^2}+ \frac{L^2_{gyy}C^2_{fx}}{\mu^4}$ and $F^* = \inf_{x\in \mathbb{R}^d}F(x) > -\infty$.
\end{theorem}

\begin{proof}

According to the line 3 of Algorithm~\ref{alg:1}, we have $u_t=\nabla_y g(x_t,y_t)$.
Then by using Lemma \ref{lem:B3}, we have
\begin{align} \label{eq:DD1}
& g(x_{t+1},y_{t+1}) - G(x_{t+1}) - \big(g(x_t,y_t) -G(x_t)\big) \nonumber \\
& \leq  -\frac{\lambda\mu}{2} \big(g(x_t,y_t) -G(x_t)\big) + \frac{1}{8\gamma}\|x_{t+1}-x_t\|^2  -\frac{1}{4\lambda }\|y_{t+1}-y_t\|^2 + \lambda\|\nabla_y g(x_t,y_t)-u_t\|^2 \nonumber \\
& = -\frac{\lambda\mu}{2} \big(g(x_t,y_t) -G(x_t)\big) + \frac{\gamma}{8}\|\widetilde{\nabla}f(x_t, y_t, v_t)\|^2
-\frac{1}{4\lambda }\|y_{t+1}-y_t\|^2,
\end{align}
where the above equality holds by $u_t=\nabla_y g(x_t,y_t)$ and $x_{t+1}=x_t-\gamma w_t=x_t-\gamma \widetilde{\nabla}f(x_t, y_t, v_t)$.

By using Lemma \ref{lem:B1}, we have
\begin{align} \label{eq:DD2}
& \|v_{t+1} - v^*_{t+1}\|^2 -\|v_t -v^*_t\|^2 \nonumber \\
& \leq  -\frac{\mu\tau}{4}\|v_t -v^*_t\|^2
 - \frac{3}{4} \|v_{t+1}-v_t\|^2
 + \frac{25\tau L^2_{gyy}r^4_v\delta^2_\epsilon}{6\mu}\nonumber \\
 & \quad + \frac{20}{3}\big(\frac{L^2_f}{\mu^2}+ \frac{L^2_{gyy}C^2_{fx}}{\mu^4}\big)\big(\|x_{t+1}-x_t\|^2 + \|y_{t+1}-y_t\|^2\big) \nonumber \\
 &=  -\frac{\mu\tau}{4}\|v_t -v^*_t\|^2
 - \frac{3}{4} \|v_{t+1}-v_t\|^2
 + \frac{25\tau L^2_{gyy}r^4_v\delta^2_\epsilon}{6\mu} + \frac{20}{3}\breve{L}^2\gamma^2\|\widetilde{\nabla}f(x_t, y_t, v_t)\|^2 + \frac{20}{3}\breve{L}^2\|y_{t+1}-y_t\|^2,
\end{align}
where the above equality holds by $x_{t+1}=x_t-\gamma w_t=x_t-\gamma \widetilde{\nabla}f(x_t, y_t, v_t)$
and $\breve{L}^2=\frac{L^2_f}{\mu^2}+ \frac{L^2_{gyy}C^2_{fx}}{\mu^4}$.
Then by using Lemma \ref{lem:B4}, we have
\begin{align} \label{eq:DD3}
    F(x_{t+1}) - F(x_{t}) & \leq - \frac{\gamma}{2}\|\nabla F(x_t)\|^2 - \frac{\gamma}{4}\|\widetilde{\nabla}f(x_t, y_t, v_t)\|^2 + \gamma L^2_{gxy}r^4_v\delta^2_{\epsilon} \nonumber \\
   & \quad + \frac{6\gamma}{\mu}\big(L^2_f+ r^2_vL^2_{gxy}\big)\big(g(x_t,y_t)-G(x_t)\big)+ 3\gamma C^2_{gxy}\|v_t-v^*_t\|^2.
\end{align}

Next, we define a useful Lyapunov function (i.e. potential function), for any $t\geq 1$
\begin{align}
 \Omega_t = F(x_t) + g(x_t,y_t)-G(x_t) + \|v_t-v_t^*\|^2.
\end{align}
By combining the above inequalities~(\ref{eq:DD1}),~(\ref{eq:DD2}) and~(\ref{eq:DD3}), we have
\begin{align} \label{eq:DD4}
 & \Omega_{t+1} - \Omega_t \nonumber \\
 & = F(x_{t+1}) - F(x_t) + g(x_{t+1},y_{t+1})-G(x_{t+1})-\big(g(x_t,y_t)-G(x_t)\big) + \|v_{t+1}-v_{t+1}^*\|^2-\|v_t-v_t^*\|^2
 \nonumber \\
 & \leq - \frac{\gamma}{2}\|\nabla F(x_t)\|^2 - \frac{\gamma}{4}\|\widetilde{\nabla}f(x_t, y_t, v_t)\|^2 + \gamma L^2_{gxy}r^4_v\delta^2_{\epsilon}+ \frac{6\gamma}{\mu}\big(L^2_f+ r^2_vL^2_{gxy}\big)\big(g(x_t,y_t)-G(x_t)\big) \nonumber \\
   & \quad + 3\gamma C^2_{gxy}\|v_t-v^*_t\|^2 -\frac{\lambda\mu}{2} \big(g(x_t,y_t) -G(x_t)\big) + \frac{\gamma}{8}\|\widetilde{\nabla}f(x_t, y_t, v_t)\|^2
-\frac{1}{4\lambda }\|y_{t+1}-y_t\|^2 \nonumber \\
& \quad -\frac{\mu\tau}{4}\|v_t -v^*_t\|^2
 - \frac{3}{4} \|v_{t+1}-v_t\|^2
 + \frac{25\tau L^2_{gyy}r^4_v\delta^2_\epsilon}{6\mu} + \frac{20}{3}\breve{L}^2\gamma^2\|\widetilde{\nabla}f(x_t, y_t, v_t)\|^2 + \frac{20}{3}\breve{L}^2\|y_{t+1}-y_t\|^2 \nonumber \\
& \leq - \frac{\gamma}{2}\|\nabla F(x_t)\|^2 -\big( \frac{\gamma}{8}-\frac{20}{3}\breve{L}^2\gamma^2\big)\|\widetilde{\nabla}f(x_t, y_t, v_t)\|^2 -\Big(\frac{\lambda\mu}{2}-\frac{6\gamma}{\mu}\big(L^2_f+ r^2_vL^2_{gxy}\big)\Big)\big(g(x_t,y_t)-G(x_t)\big) \nonumber \\
& \quad -\big(\frac{1}{4\lambda }-\frac{20}{3}\breve{L}^2\big)\|y_{t+1}-y_t\|^2 -\big(\frac{\mu\tau}{4}-3\gamma C^2_{gxy}\big)\|v_t -v^*_t\|^2 + \gamma L^2_{gxy}r^4_v\delta^2_{\epsilon}+ \frac{25\tau L^2_{gyy}r^4_v\delta^2_\epsilon}{6\mu} \nonumber \\
& \leq - \frac{\gamma}{2}\|\nabla F(x_t)\|^2 + \gamma L^2_{gxy}r^4_v\delta^2_{\epsilon}+ \frac{25\tau L^2_{gyy}r^4_v\delta^2_\epsilon}{6\mu},
\end{align}
where the last inequality holds by $0<\gamma\leq \min\Big(\frac{3}{160\breve{L}^2}, \frac{\mu\tau}{12C^2_{gxy}},\frac{\lambda\mu^2}{12\big(L^2_f+ r^2_vL^2_{gxy}\big)}\Big)$ and $0< \lambda\leq \frac{3}{80\breve{L}^2}$.
Then we can obtain
\begin{align} \label{eq:DD5}
\|\nabla F(x_t)\|^2 \leq \frac{2\big(\Omega_t  -\Omega_{t+1}\big)}{\gamma} + 2L^2_{gxy}r^4_v\delta^2_{\epsilon}+ \frac{25\tau L^2_{gyy}r^4_v\delta^2_\epsilon}{2\gamma\mu}.
\end{align}

Averaging the above inequality~(\ref{eq:DD5}), we have
\begin{align}
\frac{1}{T}\sum_{t=1}^T\|\nabla F(x_t)\|^2 & \leq \frac{2\big(\Omega_1  -\Omega_{T+1}\big)}{T\gamma} + 2L^2_{gxy}r^4_v\delta^2_{\epsilon}+ \frac{25\tau L^2_{gyy}r^4_v\delta^2_\epsilon}{2\gamma\mu} \nonumber \\
& \mathop{\leq}^{(i)} \frac{2(\Omega_1  - F^*)}{T\gamma} + 2L^2_{gxy}r^4_v\delta^2_{\epsilon}+ \frac{25\tau L^2_{gyy}r^4_v\delta^2_\epsilon}{2\gamma\mu} \nonumber \\
& = \frac{2\big(F(x_1) + g(x_1,y_1)-G(x_1) + \|v_1-v_1^*\|^2  - F^*\big)}{T\gamma} + 2L^2_{gxy}r^4_v\delta^2_{\epsilon}+ \frac{25\tau L^2_{gyy}r^4_v\delta^2_\epsilon}{2\gamma\mu},
\end{align}
where the above inequality (i) holds by $F^* = \inf_{x\in \mathbb{R}^d}F(x) > -\infty$.

Set $\delta_{\epsilon}=O\big(\frac{1}{\sqrt{T}\max(L^2_{gxy},L^2_{gyy}/\mu)r^2_v}\big)$, we can obtain
\begin{align}
\min_{1\leq t\leq T}\|\nabla F(x_t)\|^2  \leq \frac{1}{T}\sum_{t=1}^T\|\nabla F(x_t)\|^2 & \leq O(\frac{1}{T}).
\end{align}

\end{proof}

\section{Related Works}
\label{appendix-B}
The GALET method~\citep{xiao2023generalized} is meaningless for nonconvex-PL bilevel optimization, which is based on the following facts:
\begin{itemize}
\item[1)]	In the convergence analysis, the GALET method simultaneously uses the PL condition, its Assumption 2 (i.e., let $\sigma_g = \inf_{x,y}\{\sigma_{\min}^{+}(\nabla^2_{yy} g(x,y))\} >0$ for all $(x,y)$) and its Assumption 1 (i.e., $\nabla^2_{yy} g(x,y)$ is Lipschitz continuous).
\item[2)]	In the nonconvex-PL bilevel optimization problems, Hessian matrix  $\nabla^2_{yy} g(x,y))$ has two cases:  \textbf{the first case}: $\nabla^2_{yy} g(x,y)$ is singular; \textbf{the second case}: $\nabla^2_{yy} g(x,y)$ is not singular.
\item[3)] \textbf{The first case}: $\nabla^2_{yy} g(x,y)$ is singular. Since $\nabla^2_{yy} g(x,y)$ is Lipschitz continuous by Assumption 1 of~\citep{xiao2023generalized}, the singular-value of $\nabla^2_{yy} g(x,y)$ also is continuous. Thus, combining its Assumption 1 with Assumption 2 imply that the lower bound of the non-zero singular values $\sigma_g = \inf_{x,y}\{\sigma_{\min}^{+}(\nabla^2_{yy} g(x,y))\} >0$ is close to zero. Under this case, the constant $L_w = \frac{\ell_{f,1}}{\sigma_g}+\frac{\sqrt{2}\ell_{g,2}\ell_{f,0}}{\sigma_g^2}\rightarrow + \infty$ used in its Lemmas 6 and 9, and $L_F = \ell_{f,0}(\ell_{f,1}+\ell_{g,2})/\sigma_g \rightarrow + \infty$ used in its Lemma 12.
\item[4)] \textbf{The second case}: $\nabla^2_{yy} g(x,y)$ is not singular. By Assumption 2 of ~\citep{xiao2023generalized}, the singular values of Hessian is bounded away from 0, i.e., $\sigma_g>0$. Under the this case, the PL condition, Lipschitz continuous of Hessian and its Assumption 2 (the singular values of Hessian is bounded away from 0, i.e., $\sigma_g>0$) imply that GALET uses strongly convex assumption on $g(x,y)$ at variable $y$. \textbf{Note that} although the singular values of Hessian $\nabla^2_{yy} g(x,y)$ exclude zero, i.e, the eigenvalues of Hessian $\nabla^2_{yy} g(x,y)$ may be in $[-\ell_{g,2},-\sigma_g]\bigcup [\sigma_g,\ell_{g,2}]$, we cannot have negative eigenvalues at the minimizer $y^*(x)$. Meanwhile, since Hessian is Lipschitz continuous, its all eigenvalues are in $[\sigma_g,\ell_{g,2}]$. Thus, the PL condition, Lipschitz continuous of Hessian and its Assumption 2 (the singular values of Hessian is bounded away from 0, i.e., $\sigma_g>0$) imply that the GALET assumes the strongly convex.
\end{itemize}

Although \cite{kwon2023penalty} studied the nonconvex-PL bilevel optimization, it also requires some relatively strict assumptions (e.g., Assumption 1, 4,5,6,7,8 of \cite{kwon2023penalty}). For example, Assumption 1 of \cite{kwon2023penalty} gives proximal error-bound (EB) condition that is analogous to PL condition, and its Assumption 4 (2) requires the bounded $|f(x,y)|$. In particular, its Assumption 7 (2) assumes the upper-level function $ f(x,y)$ has Lipschitz Hessian. Under these conditions, \cite{kwon2023penalty} has a gradient complexity of $O(\epsilon^{-1.5})$ for finding an $ \epsilon$-stationary solution of nonconvex-PL bilevel optimization. However, without relying on Lipschitz Hessian of function $f(x,y)$ and bounded $|f(x,y)|$, our algorithm obtains an optimal gradient complexity of $O(\epsilon^{-1})$, which matches the lower bound established by the first-order method for finding an $\epsilon$-stationary point of nonconvex smooth optimization~\citep{carmon2020lower}.

Meanwhile, \cite{chen2024bilevel} studied the nonconvex-PL bilevel optimization, but it also relies on some strict assumptions, e.g.,  $h_{ \sigma}= \sigma f+ g$ is $\mu$-PL (Please see the Assumption 4.1 (a) of \cite{chen2024bilevel}). While our paper only assumes the lower-level function $g$ is $\mu$-PL. When $\sigma>0$, Assumption 4.1 (a) of \cite{chen2024bilevel} is stricter than our assumption (the lower-level function $g$ is $\mu$-PL).
In particular, the upper-level function  $ f(x,y)$ also requires Lipschitz Hessian (Please see the Assumption 4.1 (d) of \cite{chen2024bilevel}).

\textbf{Note that} from \cite{carmon2020lower}, the optimal gradient complexity is $O(\epsilon^{-0.75})$ (or $O(\epsilon^{-1.5})$) for finding an $\epsilon$-stationary point of smooth nonconvex optimization problem $\min_x f(x)$ with Lipschitz Hessian condition, i.e,  $\|\nabla f(x) \|^2\leq \epsilon $ (or $\|\nabla f(x) \|\leq \epsilon $). Meanwhile, based on Lipschitz Hessian of function $f(x,y)$, \cite{yang2023accelerating} can obtain a lower gradient complexity of $O(\epsilon^{-0.875})$ (or $O(\epsilon^{-1.75})$) for finding an $\epsilon$-stationary point of
nonconvex strongly-convex bilevel optimization, i.e., $\|\nabla \Phi(x) \|^2\leq \epsilon$ (or $\|\nabla \Phi(x) \|\leq \epsilon $).
Under these strict assumptions, thus, although \cite{chen2024bilevel} also obtain a gradient complexity of $O(\epsilon^{-1})$,
 \emph{this is not optimal gradient complexity.}

\begin{lemma}
 (Lemma G.6 of \cite{chen2024bilevel})
 For a $\mu$-PL function $h(x):\mathbb{R}^d \rightarrow \mathbb{R}$  that is twice differentiable, at any $x^* \in \arg\min_x h(x)$,
\begin{align}
 \lambda^+_{\min} (\nabla^2 h(x^*)) \geq \mu,
\end{align}
where $\lambda^+_{\min} (\cdot)$ denotes the smallest non-zero eigenvalue.
\end{lemma}

 In fact, Lemma G.6 of \cite{chen2024bilevel} is exactly useful for our HJFBiO method. Based on Lemma G.6 of \cite{chen2024bilevel}, our Assumption~\ref{ass:2} is reasonable when
 has an unique minimizer. Meanwhile, our Assumption~\ref{ass:2g} also is reasonable when
 have multiple local minimizers.

\end{document}